\documentclass[a4paper,12pt]{article}
\usepackage{amsmath,amsfonts,amssymb,amsthm}
\usepackage{graphics}
\usepackage{epsfig}
\usepackage{color}
\definecolor{darkred}{rgb}{0.9,0.1,0.1}

\newcommand{\ignore}[1]{}
\newcommand{\av}{-\hspace{-2.4ex}\int}
\newcommand{\an}{\alpha_0}
\newcommand{\Dprime}{{D'}^2(u,1)}

\newcommand{\te}{\tilde{\eta}}
\newcommand{\px}{\partial_x}

\newtheorem{proposition}{Proposition}
\newtheorem{theorem}{Theorem}

\newtheorem{lemma}{Lemma}
\newtheorem{corollary}{Corollary}
\usepackage{stmaryrd}  

\begin{document}

\begin{center}
\begin{large}
H\"older regularity for a non-linear parabolic equation
driven by space-time white noise\\
F.\ Otto, H.\ Weber\\
\end{large}
\end{center}

\begin{center}
\begin{minipage}{13cm}
\noindent {\bf Abstract.}
We consider the non-linear equation $T^{-1} u+\partial_tu-\partial_x^2\pi(u)=\xi$
driven by space-time white noise $\xi$, which is uniformly parabolic because 
we assume that $\pi'$ is bounded away from zero and infinity. Under the further assumption
of Lipschitz continuity of $\pi'$ we show that the stationary solution is --- as for the linear case ---
almost surely H\"older continuous with exponent
$\alpha$ for any $\alpha<\frac{1}{2}$ w.\ r.\ t.\ the parabolic metric. More precisely,
we show that the corresponding local H\"older norm has stretched exponential moments.

\smallskip

On the stochastic side, we use a combination of 
martingale arguments to get second moment estimates with concentration of measure
arguments to upgrade to Gaussian moments. On the deterministic side, we first perform 
a Campanato iteration based on the De Giorgi-Nash Theorem as well as finite and 
infinitesimal versions of the $H^{-1}$-contraction principle, which yields Gaussian moments 
for a weaker H\"older norm. In a second step this estimate is improved to the optimal 
H\"older exponent at the expense of weakening the integrability to stretched exponential.

\end{minipage}
\end{center}

\section{Introduction and main result}

We are interested in the stochastic nonlinear parabolic equation 
\begin{equation}\label{i.1}
T^{-1}u+\partial_tu-\partial_x^2\pi(u)=\xi,
\end{equation}
where $\xi$ denotes space-time white noise. The nonlinear character of (\ref{i.1}) is that of
a fully nonlinear equation rather than a quasi-linear equation, since rewriting (\ref{i.1}) as
the quasi-linear equation (\ref{i.1bis}) is not helpful as we explain below, and since the 
deterministic estimates we need are related to the linearization of a fully nonlinear equation,
cf.\ (\ref{i.5}), rather than to the linearization of a quasi-linear equation
(this distinction would be more pronounced in a multi-dimensional case).
We assume that the nonlinearity $\pi$ is uniformly
elliptic in the sense that there exists a $\lambda>0$ such that
\begin{equation}\label{i.2}
\lambda\le \pi'(u)\le 1\quad\mbox{for all}\;u\in\mathbb{R}.
\end{equation}
In particular, this rules out the degenerate case that goes under the name of porous medium equation.
Furthermore, we assume some regularity of $\pi$ in the sense that there exists $L<\infty$
such that
\begin{equation}\label{i.3}
|\pi''(u)|\le L\quad\mbox{for all}\;u\in\mathbb{R}.
\end{equation}
We are interested in H\"older regularity of solutions of (\ref{i.1}); the simplest solution
to (\ref{i.1}) is the space-time stationary solution $u$ of (\ref{i.1}) on which we shall focus
in this paper. The main reason for including 
the massive term in (\ref{i.1}) (i.\ e.\ assuming $T<\infty$) is to ensure existence and uniqueness of this object; 
the only other role is to provide a large-scale estimate through Lemma \ref{L0}. 
In this version of the paper, we will be {\it completely informal} about why and in which sense (\ref{i.1})
is well-posed, and why the martingale and concentration of measure arguments can be carried out
(we will just motivate them when we first need them).

\medskip

A crucial insight is that the law of the (unique) stationary random field
$u$ is invariant under the rescaling
\begin{equation}\label{t1}
x=R\hat x,\quad t=R^2\hat t,\quad u=R^\frac{1}{2}\hat u,
\end{equation}
provided, the nonlinearity and the massive term are adjusted according to
\begin{equation}\label{t2}
\hat\pi(\hat u)=R^{-\frac{1}{2}}\pi(R^\frac{1}{2}\hat u),
\quad \hat T=R^{-2}T.
\end{equation}
For this observation we used that in view of its defining relation
$$\big\langle \big(\int \zeta \xi dxdt  \big)^2 \big\rangle = \int \zeta^2 dx dt$$ for a test function $\zeta$
(that is, loosely speaking $\int\zeta\langle\xi(t,x)\xi(0,0)\rangle dxdt =\zeta(0,0)$), 
space-time white noise rescales as $\xi=\frac{1}{\sqrt{R R^2}}\hat\xi=R^{-\frac{3}{2}}\hat\xi$.
From this invariance property we learn that as we go to small scales (i.\ e.\ $R\ll 1$),
the effective nonlinearity as measured by the Lipschitz constant $L$ of $\pi'$ in
(\ref{i.3}) decreases according to 
\begin{equation}\label{i.b}
\hat L=R^\frac{1}{2}L.
\end{equation}
This suggests that
on small scales, $u$ has the same regularity as if (\ref{i.1}) were replaced by
its linear version (without massive term) $\partial_t u-a_0\partial_x^2u=\xi$ for some
constant $a_0\in[\lambda,1]$. Hence
we expect that on small scales, $u$ is H\"older continuous with exponents $\alpha$
(in the parabolic Carnot-Carath\'eodory geometry) for any $\alpha<\frac{1}{2}$. This is
exactly what we show, making crucial use of the above scale invariance. 

\medskip

We note in passing that it is {\it not} helpful to write the elliptic operator in the more 
symmetric form
\begin{equation}\label{i.1bis}
T^{-1}u+\partial_tu-\partial_x(\pi'(u)\partial_xu)=\xi,
\end{equation}
since even in case of the stochastic heat equation, $u$ (and thus $\pi'(u)$) is a function
in the H\"older space with exponent $\frac{1}{2}-$ so that $\partial_xu$ would be a distribution
in the (negative) H\"older space with exponent $-\frac{1}{2}-$, so that there is no standard
distributional definition of the product $\pi'(u)\partial_xu$. In fact,
rather than appealing to regularity theory for linear but non-constant
coefficient equations of the form 
$T^{-1}u+\partial_tu-\partial_x(a\partial_xu)=\partial_xg$,
we have to appeal to the theory for
$T^{-1}w+\partial_tw-\partial_x^2(a w)=\partial_x^2g$,
cf.\ Proposition \ref{P4}. 

\medskip

Let us now briefly comment on existing regularity theory for non-linear parabolic stochastic differential
equation of the type of (\ref{i.1}).
There is a large body of literature on stochastic equations of the type (\ref{i.1}), but mostly with a quite
different focus: The focus there is to tackle on the one hand more nonlinear situations, 
like the case of a degenerate ellipticity
(i.\ e.\ $\lambda=0$ in (\ref{i.2})) or the case of multiplicative noise, but on the other hand
to assume ``whatever it takes'' on the spatial covariance structure of the noise. Sometimes, structural assumptions allow to 
mimic an approach that is obvious in the semi-linear case, namely the approach of
decomposing the solution into a rough part $w$ that solves a more explicitly treatable
stochastic differential equation and a more regular part $v$ that 
solves a parabolic equation with random coefficients and/or right-hand-side described through $w$,
and then allows for an application of deterministic regularity theory. We refer to
\cite{Gess} for an example with a multiplicative decomposition of this type.
The recent work by Debussche et.\ al.\ on quasi-linear parabolic stochastic equations, i.\ e.\
equations of the form (\ref{i.1bis}) or more generally with an elliptic operator of the form $-\nabla\cdot a(u)\nabla u$, 
refines this approach to a fixed point argument,
and appeals to the De Giorgi-Nash Theorem, which yields a H\"older a priori bound on linear parabolic equations with just uniformly elliptic
coefficients as a starting point to bootstrap to the optimal H\"older continuity via Schauder theory,
see \cite[Introduction]{Debussche}. However, cf.\ the above discussion of
(\ref{i.1bis}), this treatment seems limited to situations where the noise $\xi$
is so regular that in the case of the linear equation, $\nabla u$ is at least locally integrable in time-space
(to be more quantitative: $\int|\nabla u|^pdxdt<\infty$ for some $p>3$ on the level of space-time
isotropic $L^p$-norms).

\medskip

By the equivalence of Campanato and H\"older spaces, see for instance \cite[Theorem 5.5]{Giaquinta}, 
H\"older continuity can be expressed in terms of a localized $L^2$-modulus of continuity.
Because of the eventual conditioning on the distant noise, it is more convenient to replace
a sharp spatial localization on parabolic cylinders
by a soft localization via an exponentially decaying function
\begin{equation}\label{p31}
\eta(x)=\frac{1}{2}\exp(-|x|),\quad\eta_r(x):=\frac{1}{r}\eta(\frac{x}{r}),
\end{equation}
note that the normalization imply that $\int\eta_r\cdot dx$ corresponds to a spatial average
that is localized near the origin on scale $r$. We note that while the exponential form
of the cut-off is probably not essential (any thicker than Gaussian tails should suffice),
it is convenient at many places of the proof.
Abbreviating the $L^2$-modulus of continuity at the origin and on parabolic scale $r$ by
\begin{equation}\nonumber
D^2(u,r):=\av_{-r^2}^0\int\eta_r(u-\av_{-r^2}^0\int\eta_ru)^2dxdt,
\end{equation}
our first main result reads as follows:

\begin{theorem}\label{T}
Let $u$ be the unique  stationary solution to \eqref{i.1}.
W.\ l.\ o.\ g.\ suppose that $T=1$ in (\ref{i.1}). There exists a H\"older exponent $\an \in(0,\frac{1}{2})$ 
depending only on $\lambda$, so that we have a Gaussian bound for the $\alpha_0$-H\"older $L^2$-averaged
modulus of continuity at the origin in the sense that
\begin{equation}\label{t7}
\Big\langle\exp\Big(\frac{1}{C}\big(\sup_{r\le 1}\frac{1}{r^{\alpha_0} }D(u,r)\big)^{2}
\Big)\Big\rangle\le 2,
\end{equation}
with a constant $C<\infty$ only depending on  $\lambda>0$, $L<\infty$ and $\an$.
\end{theorem}
The H\"older exponent $\an\in (0,\frac12)$ is determined by an application of the celebrated De Giorgi-Nash Theorem via Proposition~\ref{P1} below (in 
fact the $\alpha_0$ in Theorem~\ref{T} is slightly smaller than the $\alpha_0$ in Proposition~\ref{P1}).   It
only depends on the ellipticity ratio $\lambda$. In our second main result we improve the H\"older regularity exponent
up to the optimal value $\alpha =\frac12-$ at the expense of weakening the integrability.  

\begin{theorem}\label{T2}
Let $u$ be the unique stationary solution to \eqref{i.1} for $T =1$. Let $\an \in (0,  \frac12)$ be 
the H\"older exponent appearing in Theorem~\ref{T}. Then for any H\"older exponent $\alpha\in(\alpha_0,\frac{1}{2})$ 
we get 
\begin{equation}\label{t7new}
\Big\langle\exp\Big( \frac{1}{C}\Big(\sup_{r\le 1}\frac{1}{r^\alpha}D(u,r)\Big)^{2 \frac{\an}{\alpha} }
\Big)\Big\rangle\le 2,
\end{equation}
with a constant $C<\infty$ only depending on $\lambda>0$, $L<\infty$  and  $\alpha<\frac{1}{2}$.
\end{theorem}

Theorems~\ref{T} and \ref{T2} imply  bounds 
for the more conventional local H\"older semi-norms of the random field $u$. 
For any $\alpha \in (0,1)$ we set
\begin{align}\label{DefHol}
[u]_{\alpha} = \sup_{R \in (0,1)} \frac{1}{R^\alpha} \sup_{\substack{(t,x), (s,y) \in (-1,0) \times (-1,1)  \\ \sqrt{|t-s|} +|x-y| <R }} |u(t,x) - u(s,y)| .
\end{align}
Theorem~\ref{T2} implies the following:
\begin{corollary}\label{C1}
Under the assumptions of Theorem~\ref{T2} we have
\begin{equation*}
\Big\langle \exp \Big( \frac{1}{C}[u]_{\alpha}^{2 \frac{\an}{\alpha} }   \Big) \Big\rangle \leq 2 \,
\end{equation*}
 for a constant $C<\infty$ which only depends on $\lambda>0$,  $L<\infty$,   $\alpha < \frac12$ and $\epsilon>0$.
\end{corollary}

\section{Strategy of proof and ingredients}\label{S2}

Theorem \ref{T}, like Lemma \ref{L0} below,
relies on a concentration of measure argument for Lipschitz random variables: 
For any a random variable $F$ that is 1-Lipschitz when considered as a path-wise functional
of the white noise $\xi$, one has $\langle\exp(\lambda F)\rangle\le\exp(\lambda\langle F\rangle+\frac{1}{2}\lambda^2)$
for any number $\lambda$. In particular, if $F\ge 0$ is 1-Lipschitz and satisfies $\langle F\rangle\le 1$,
it has Gaussian moments $\langle\exp(\frac{1}{C}F^2)\rangle\le 2$, for some universal constant $C<\infty$.
Here the norm underlying the Lipschitz property is the norm of the Cameron-Martin space, which simply
means that infinitesimal variations $\delta\xi$ of the space-time white noise
are measured in the space-time $L^2$-norm. To continue with the name-dropping,
this type of Lipschitz continuity means that the carr\'e-du-champs $|\nabla F|^2$
of the Malliavin derivative is bounded independently of the given realization
of the noise, where for a given realization $\xi$ of the noise,
$|\nabla F|$ is the smallest constant $\Lambda$ in
\begin{equation}\label{i.a}
|\delta F|\le \Lambda\Big(\int(\delta\xi)^2dxdt\Big)^\frac{1}{2}.
\end{equation}
Here $\delta F$ denotes the infinitesimal variation of $F$ generated by the infinitesimal
variation $\delta\xi$ of the noise $\xi$, a linear relation captured by the Fr\'echet derivative 
(a linear form) of $F$ w.\ r.\ t.\ $\xi$. For those not confident in this continuum version of
concentration of measure we derive it from the discrete case in the proof of Lemma \ref{L0},
where also the type of martingale arguments entering Proposition \ref{P1} via
Lemma \ref{L1} (and Lemmas \ref{L2} and \ref{L4} again) is explained for the non-expert.

\smallskip

Both in the proofs of Theorem~\ref{T} and Theorem~\ref{T2}, concentration of measure will be applied to the random variable $F=D(u,r)$.
It is Proposition \ref{P2} which provides the bound on the Malliavin derivative w.\ r.\ t.\ to
the ensemble $\langle\cdot\rangle_1$ that describes the space time white noise $\xi$
{\it restricted} to the time slice $(t,x)\in(-1,0)\times\mathbb{R}$. In particular, this
means that the admissible variations
$\delta\xi$ in (\ref{i.a}) are supported in $(t,x)\in(-1,0)\times\mathbb{R}$; we denote by $|\nabla F|_1^2$
the corresponding carr\'e-du-champs. Proposition \ref{P1} in turn provides the estimate
of the (conditional) expectation, that is, the expectation in $\langle\cdot\rangle_1$ which is used in Theorem~\ref{T}, 
while the proof of Theorem~\ref{T2} relies on Proposition~\ref{P1old}.

\smallskip

Proposition~\ref{P1} provides a bound on the expectation  of $\frac{1}{r^{\alpha_0}}D(u,r)$ in terms of quantities 
that are \emph{linear} in $D'(u,1)$ (which roughly behaves as  
$D(u,1)$) and in combination with Proposition~\ref{P2} this bound can be upgraded to Gaussian moments. 
Using the scale invariance  \eqref{t1} and \eqref{t2} this estimate can then  be used in a (stochastic) Campanato 
iteration leading to Theorem~\ref{T}. The drawback of Proposition~\ref{P1} is that it is restricted to the 
small H\"older exponent $\alpha_0$ because the proof relies on the De Giorgi-Nash Theorem. 
Proposition \ref{P1old} (in conjunction with Proposition \ref{P2}) in turn yields Gaussian moments
for $\frac{1}{r^\alpha}D(u,r)$, for all admissible exponents $\alpha \in (0, \frac12)$,
however only up to $1+r^\frac{1}{2}D'(u,1)$, which roughly behaves as 
$1+r^\frac{1}{2}D(u,1)$, and modulo the multiplicative (and nonlinear) error of $\frac{L}{r^\frac{3}{2}}(1+D(u,1))$.
Evoking the scale invariance (\ref{t1}) \& (\ref{t2}), this estimate will be used for small scales, where 
thanks to the behavior (\ref{i.b}) of $L$, the multiplicative error fades away, so that
Theorem~\ref{T2} can be obtained by another (non-linear) Campanato iteration. Theorem~\ref{T}
is needed as an anchoring for this iteration.

\medskip
%
%
%
Since these propositions will be applied to small scales,
so that in view of (\ref{t2}) the massive term fades away, we cannot expect
help from it; as a matter of fact, we will ignore the massive term in the proof (besides in
Lemma \ref{L0} where it is essential).

\begin{proposition}\label{P1}
There exists a H\"older exponent $\an\in(0,\frac{1}{2})$, depending only on $\lambda$, such 
that we have all $r\le 1$
\begin{equation}\nonumber
\langle D(u,r)\rangle_1\lesssim r^{\an} \Big( 1+ D'(u,1)  +\frac{L}{r^\frac{3}{2}}(D'(u,1)+1) \Big),
\end{equation}
where $D'(u,1)$ depends only on $u(t=-1,\cdot)$:
\begin{equation}\label{Dprime}
\Dprime:=\int\eta(u-\int\eta u)^2dx_{|t=-1}.
\end{equation}
Here and in the proof, $\lesssim$ means up to a constant only depending on $\lambda>0$.
\end{proposition}

\begin{proposition}\label{P2}
We have for the carr\'e-du-champs of the Malliavin 
derivative
\begin{equation}\nonumber
|\nabla D(u,r)|_1\lesssim r^{-\frac{3}{2}},
\end{equation}
where here and in the proof,
$\lesssim$ means up to a constant only depending on $\lambda>0$.
\end{proposition}

The only purpose of the presence of the massive term is that in the original
scale, it provides control of the $L^2$-averaged H\"older continuity on scales 1,
and thus the anchoring for the Campanato iteration:

\begin{lemma}\label{L0}
Suppose that $T=1$ in (\ref{i.1}). Then we have
\begin{equation}\nonumber
\big\langle\exp\big(\frac{1}{C} D^2(u,1)\big)\big\rangle\le 2
\end{equation}
for some constant $C$ only depending on $\lambda$.
\end{lemma}

\medskip

In order to derive Propositions \ref{P1} and \ref{P2}, we will consider differences of
solutions to (\ref{i.1}) for Proposition \ref{P1}, or infinitesimal perturbations
of solutions for Proposition \ref{P2}. Finite or infinitesimal differences of solutions
satisfy a formally linear parabolic equation with an inhomogeneous coefficient field $a$,
which in view of (\ref{i.2}) is uniformly elliptic:
\begin{equation}\label{p21}
\lambda\le a(t,x)\le 1\quad\mbox{for all}\;(t,x)\in(-1,0)\times\mathbb{R}.
\end{equation}
The linearized operator comes in the conservative form of $\partial_tu-\partial_x^2(au)$.
For a priori estimates of the corresponding initial value problem, it is most
natural to write the r.\ h.\ s.\ also in conservative form:
\begin{equation}\label{i.5}
\partial_tw-\partial_x^2(aw)=\partial_th+\partial_x^2g.
\end{equation}
The $L^2$-estimates on solutions of (\ref{i.5}) from Proposition \ref{P4}
%
%
might be seen as an infinitesimal version of the $\dot H^{-1}$-contraction principle for the
deterministic counterpart of (\ref{i.1}), which will be explicitly used in Lemma \ref{L1},
see the proof of Lemma \ref{L0}, which is a good starting point for the PDE arguments, too.

\medskip

\begin{proposition}\label{P4}
Consider a solution $w$ of (\ref{i.5}) with r.\ h.\ s.\ described by $(g,h)$. In the case
$h=0$ we get the local estimate
\begin{align}
\notag
&\sup_{t \in (-1,0)}  \int \sqrt{\eta}w  (1-\px^2)^{-1} \sqrt{\eta}w  dx+  \int_{-1}^0\int\eta w^2dxdt \\
\label{p43A}
& \qquad \lesssim   \int\sqrt{\eta}w  (1-\px^2)^{-1} \sqrt{\eta}w  dx|_{t=-1} + \int_{-1}^0\int\eta g^2dxdt. 
\end{align}
In the case of general $h$ and homogeneous initial data i.e. if 
\begin{align*}
 w=h=0 \qquad 
 \text{for } \quad t=-1 
\end{align*}
we get both the local bound
\begin{align}
& \int_{-1}^0\int\eta w^2dxdt 
\label{p43}
\lesssim   \int_{-1}^0\int\eta (g^2+h^2)dxdt
\end{align}
and the global bound
\begin{equation}\label{p24}
\int_{-1}^0\int w^2dxdt\lesssim\int_{-1}^0\int(g^2+h^2)dxdt.
\end{equation}
Here and in the proof $\ll$ and $\lesssim$ refer just to $\lambda$.
\end{proposition}
However, next to this ``soft'' a priori estimate for solutions of the initial value problem
for (\ref{i.5}), we also need the following ``hard'' a priori estimate. Well beyond the $L^2$-bound
in Proposition \ref{P4}, Proposition \ref{P5} provides equi-integrability of $v^2$ for a solution
of the homogeneous version of (\ref{i.5}) in the sense of a Morrey norm, for the latter see \cite[Definition 5.1]{Giaquinta}. 
It does so in a quite quantified way: (\ref{i.10})
provides equi-integrability as if $v\in L^p$ (in time-space) for $p=\frac{3}{1-\an}>3$.  
Loosely speaking, this equi-integrability arises as follows: The spatial anti-derivative $V$
of the given solution $v$ of the homogeneous version of (\ref{i.5}) satisfies the {\it divergence-form}
equation $\partial_tV-\partial_x(a\partial_xV)=0$. Now the celebrated theory of De Giorgi and Nash,
in particular in the parabolic version of Nash \cite{Nash}, implies that for some $\an>0$
depending only on the ellipticity $\lambda$ in (\ref{p21}), the $\an$-H\"older norm
of $V$ is controlled by weaker norms of $V$. The connection between Nash's result and (\ref{i.10})
is obvious on the level of scaling: The $\an$-H\"older norm of $V$ has the same (parabolic) scaling
as the $L^p$-norm of $v=\partial_xV$ with $\an=1-\frac{3}{p}$. On the level of equi-integrability
of $(\partial_xV)^2$, 
this scaling analogy indeed can be made rigorous with help of the parabolic Caccioppoli estimate for $V$.

\begin{proposition}\label{P5} 
Let $v$ be a solution to \eqref{i.5} with right hand side $h=f=0$. Then there exists an $\an >0$, depending only on $\lambda$ such that  for $0<r\leq1$ we have the estimate
\begin{equation}
\av_{r^2}^0 \int \eta_r v^2 \, dx dt \lesssim r^{-2+2\an} \int  \sqrt{\eta}  v  (1-\px^2)^{-1} \sqrt{\eta} v  \, dx \big|_{t=-1}.  \label{i.10}
\end{equation}
Here and in the proof $\lesssim$ only refers to $\lambda$.
\end{proposition}

The following proposition provides the non-linear version of Proposition~\ref{P1} which is used in the proof of 
Theorem~\ref{T2}. In this proposition the restriction on the H\"older exponents is removed and all exponents 
$\alpha < \frac12$ are admitted.
\begin{proposition}\label{P1old}
Pick a H\"older exponent $\alpha\in(0,\frac{1}{2})$. 
Then we have all $r\le 1$
\begin{equation}\nonumber
\langle D(u,r)\rangle_1\lesssim r^\alpha\big(1+\frac{L}{r^\frac{3}{2}}(D'(u,1)+1)\big)
\big(1+r^\frac{1}{2}D'(u,1)\big),
\end{equation}
where $D'(u,1)$ is defined in \eqref{Dprime}. 
Here and in the proof, $\lesssim$ means up to a constant only depending on $\lambda>0$ and $\alpha<\frac{1}{2}$.
\end{proposition}

For this proposition we follow a standard approach in Schauder theory for parabolic (and elliptic) equations 
and consider (\ref{i.5}) with constant coefficients $a_0\in[\lambda,1]$,
which will arise from locally ``freezing'' the variable coefficient field $a$:
\begin{equation}\label{i.6}
\partial_tv-a_0\partial_x^2v=f.
\end{equation}
Proposition \ref{P3} states a classical $L^\infty$ estimate
for (\ref{i.6}), the only difficulty coming from the low regularity of the initial
data $v_{|t=-1}$ and the moderate regularity of the r.\ h.\ s.\ $f$ assumed in Proposition \ref{P3}.
We give a self-contained proof.

\begin{proposition}\label{P3}
Consider a solution $v$ of (\ref{i.6}) with r.\ h.\ s.\ $f$. Then we have a
localized $L^\infty$-estimate
\begin{equation}\label{p35}
\sup_{(t,x)\in(-1,0)\times\mathbb{R}}(t+1)^\frac{1}{2}\eta v^2\lesssim
\int_{-1}^0\int\eta f^2dxdt+\int\eta v^2dx_{|t=-1}.
\end{equation}
%
\end{proposition}

\medskip

We'd like to point out a synergy in terms of methods between this approach to
regularity for stochastic partial differential equations driven by stationary noise,
and an approach to regularity for elliptic partial differential equations with stationary random coefficient field
that is emerging over the past years \cite{MarahrensOtto,ArmstrongSmart,GloriaNeukammOtto}.
At first glance, the differences dominate:
Here, we have a {\it nonlinear} and {\it parabolic} partial differential equation driven by a random
right-hand-side $\xi$, and we hope for almost-sure
{\it small-scale} regularity {\it despite the short-range decorrelation} of $\xi$,
which implies its roughness.
There, the main features already appear on the level of a {\it linear} and {\it elliptic}
equation, for instance on the level of the harmonic coordinates or the corrector $\phi_i$ given by 
$-\nabla\cdot a(\nabla \phi_i+e_i)=0$ where $e_i$ is the $i$-th unit vector, and one hopes
for almost-sure {\it large scale} regularity {\it thanks to the long-range decorrelation} of the coefficient field
$a$.
In the first case, {\it randomness limits} H\"older regularity, whereas in the second case, 
{\it randomness improves}
H\"older regularity: In fact, for almost every realization of $a$, $a$-harmonic functions $u$ satisfy
a first-order Liouville principle \cite{GloriaNeukammOtto}, and even Liouville principles
of any order \cite{FischerOtto}, which is the
simplest way to encode large-scale H\"older regularity. Even the lowest-order Liouville principle is known to fail
for some uniformly elliptic and smooth coefficient fields $a$, so that these results
indeed show a regularizing effect of randomness. 

\medskip

Despite these obvious differences, the approach is very similar: Both here and there
(in \cite{MarahrensOtto} and, more explicitly, in \cite{FischerOtto2}, \cite{GloriaNeukammOtto}) 
 one is appealing to
the combination of sensitivity estimates (how do certain functionals of the solution 
depend on the right hand side here, or on the coefficient field there?) measured in terms
of a carr\'e du champs (of the Malliavin derivative here, 
or of a suitable vertical derivative that is compatible with the correlation structure there),
and then appeals to concentration of measure (on the Gaussian level here, or via the intermediate
of a Logarithmic Sobolev Inequality there).

\medskip

Such a synergy in methods that treat models with thermal noise like in high- or infinite dimensional 
stochastic differential equations with reversible invariant (Gibbs) measure and those that treat models with
quenched noise like in stochastic homogenization is not new: In their seminal work on Gradient Gibbs measures, 
a model in statistical mechanics that describes thermally fluctuating surfaces, Naddaf and Spencer appeal 
to stochastic homogenization to characterize the large-scale correlation structure of the field \cite{NaddafSpencer}. 
Their analysis can also be interpreted as considering the infinite-dimensional stochastic differential 
equation of which the measure is the reversible invariant measure, an equation which can be seen as a 
spatial discretization of a stochastic nonlinear parabolic partial differential equation, 
and to consider the Malliavin derivative of its solution with respect to the (discrete) space-time white noise
\cite{Conlon}. Again, the nonlinearity is rather of the symmetric form (\ref{i.1bis}) 
and Naddaf and Spencer appeal to Nash's heat kernel bounds.

\medskip

We close this parenthesis by noting that for stochastic partial differential equations
and stochastic homogenization, even the deterministic ingredients are similar: In both cases, the sensitivity estimate 
leads to a {\it linear} partial differential equation 
(parabolic here, elliptic there) with a priori only
uniformly elliptic coefficient field (in space-time here, in space there), that is,
without any a priori modulus of continuity. 
In both cases, a buckling argument is needed to obtain bounds on H\"older norms with high
stochastic integrability.
While here, the need of a buckling estimate is obvious since the small-scale regularity 
of the coefficient field $a=\pi'(u)$ in the sensitivity equation
is determined by the small-scale regularity of the solution $u$ around which one is linearizing,
the buckling is less obvious there: It turns out that the large-scale regularity properties of the
operator $-\nabla\cdot a\nabla$ are determined by the large-scale properties 
of the harmonic coordinates $x_i+\phi_i$, the special solution mentioned above.
Here, buckling proceed by showing that the linear operator $\partial_t-\partial_x^2a$ is close to a 
constant coefficient operator $\partial_t-a_0\partial_x^2$
on small scales, there, it proceeds by showing that it is close to a constant coefficient operator
on large scales, namely the homogenized operator $-\nabla\cdot a_{hom}\nabla$. 
In both cases, a Campanato-type iteration is the appropriate deterministic tool for the buckling.
Here, this is not surprising since Campanato iteration is a robust way of deriving H\"older estimates
(see for instance \cite[Chapter 5]{Giaquinta});
there, the use of Campanato iteration to push the constant-coefficient regularity theory from
the infinite scale to large but finite scales was first introduced in \cite{AvellanedaLin} in case
of periodic homogenization, then transferred to stochastic homogenization in \cite{ArmstrongSmart},
and refined in \cite{GloriaNeukammOtto} in a way that brings it very close to its small-scale application. 

\medskip

After this aside, we turn back to our proof.
Next to these deterministic ingredients, Proposition \ref{P1} and Proposition~\ref{P1old} also require
a couple of classical, second moment stochastic estimates. The first lemma
provides such a low-stochastic moment estimate on the $L^2$-H\"older-$\frac{1}{2}$ 
modulus of continuity, which however is restricted to a spatial modulus
and is only localized to scales $1$.
This spatial $L^2$-H\"older modulus of continuity is expressed in terms of the $L^2$-difference
of spatial shifts (which are then exponentially averaged over the shifts); this
form arises naturally from a martingale version of the (deterministic) $\dot H^{-1}$-contraction principle
for equations of the form (\ref{i.1}) with uniform ellipticity (\ref{i.2}). 
In fact, we use a spatially localized version of the $\dot H^{-1}$-contraction
principle.

\begin{lemma}\label{L1} Let $u$ denote the stationary solution of (\ref{i.1}) and denote by
$u^h$ its spatial translation by the shift $h\in\mathbb{R}$. Then we have for $r\ll 1$
\begin{eqnarray}
\lefteqn{\Big\langle\int\eta_r(h)\int_{-\frac{1}{2}}^0\int\eta(u^h-u)^2dxdtdh\Big\rangle_1}\nonumber\\
&\lesssim& r+r^2\Big\langle\int_{-1}^0\int\eta(u-\int_{-1}^0\int\eta u)^2dxdt\Big\rangle_1
=r+r^2\langle D^2(u,1)\rangle_1.\nonumber
\end{eqnarray}
Here and in the proof $\lesssim$ and $\ll$ just refer to $\lambda$.
\end{lemma}

The second (very similar) step is to estimate the ``bulk'' $L^2$-modulus on the r.\ h.\ s.\ of Lemma \ref{L1}
by the boundary $L^2$-modulus of the initial data $u(t=-1,\cdot)$.

\begin{lemma}\label{L2} The stationary solution $u$ of (\ref{i.1}) satisfies
\begin{eqnarray}
\lefteqn{\langle D^2(u,1)\rangle_1
=\Big\langle\int_{-1}^0\int\eta(u-\int_{-1}^0\int\eta u)^2dxdt\Big\rangle_1}\nonumber\\
&\lesssim& 1+\int\eta(u-\int\eta u)^2dx_{|t=-1}=1+\Dprime.
\end{eqnarray}
Here and in the proof $\lesssim$ and $\ll$ just refer to $\lambda$.
\end{lemma}

The third step is to upgrade the purely spatial $L^2$-averaged H\"older-$\frac{1}{2}$
modulus of continuity into a space-{\it time} modulus of continuity. 

\begin{lemma}\label{L4}
The stationary solution $u$ of (\ref{i.1}) satisfies for $r\ll 1$
\begin{eqnarray*}
\lefteqn{\Big\langle\Big(\av_{-r^2}^0(\int\eta_ru
-\av_{-r^2}^0\int\eta_ru)^2dxdt\Big)^\frac{1}{2}\Big\rangle_1}\\
&\lesssim&r^\frac{1}{2}
+\Big\langle\Big(\av_{-r^2}^0\int\eta_r(u-\int\eta_ru)^2 dxdt\Big)^\frac{1}{2}\Big\rangle_1.
\end{eqnarray*}
Here and in the proof $\lesssim$ and $\ll$ just refer to $\lambda$.
\end{lemma}

The crucial ingredient for Proposition \ref{P1} is the passage from 
measuring the H\"older-$\alpha$ $L^2$-modulus of continuity
on scales 1 down to scales $r$. It is here that we need the deterministic ingredients
of Propositions \ref{P3} and \ref{P4}.
Not surprisingly, we will need in this argument that solutions $g$ to the stochastic {\it linear
constant coefficient} parabolic equation, around which we perturb, have this localization
property. This is provided by the following localized space-time {\it supremum} estimate of the 
H\"older-$\alpha$ modulus of continuity of $g$.

\begin{lemma}\label{L3} For $a_0 \in [\lambda,1]$ let $g(a_0, \cdot, \cdot)$ be the solution of 
\begin{equation}\nonumber
\partial_tg-a_0\partial_x^2g=\xi\;\;\mbox{for}\; t>-1,\quad
g(t=-1,\cdot)=0
\end{equation}
%
Then for any H\"older exponent $\alpha<\frac{1}{2}$ 
and all shifts $h\in\mathbb{R}$
\begin{equation}\nonumber
\Big\langle\sup_{a_0 \in [\lambda,1]}\sup_{(t,x) \in(-1,0)\times\mathbb{R}}\eta(g^h-g)^2\Big\rangle_1
\lesssim\min\{|h|^{2\alpha},1\}.
\end{equation}
Here and in the proof $\lesssim$ and $\ll$ refer to $\lambda$ and $\alpha$.
\end{lemma}


\section{Proofs}

We start the string of proofs with Lemma \ref{L0}, since it contains the other
arguments {\it in nuce}.

{\sc Proof of Lemma \ref{L0}}.
We will establish the lemma in the stronger version where instead of
$D^2(u,1)$, we control the Gaussian moments of $E^2(u,1):=\int_{-1}^0\int\eta u^2dxdt\ge D^2(u,1)$:
\begin{equation}\nonumber
\Big\langle\exp\big(\frac{1}{C}E^2(u)\big)\Big\rangle\lesssim1.
\end{equation}
By concentration of measure, cf.\ beginning of Section \ref{S2}, this is a consequence of the bound on the expectation
\begin{equation}\label{l04}
\langle E^2(u,1)\rangle\lesssim1
\end{equation}
and the uniform bound on the carr\'e-du-champs of the Malliavin derivative
\begin{equation}\label{l01}
|\nabla E(u,1)|^2\lesssim 1.
\end{equation}
In order gain confidence in this principle of concentration of measure,
let us relate it to the discrete case, that is, the case of countably
many independent normal Gaussian random variables, see for instance \cite[p.135]{Ledoux}
for a proof of concentration of measure by an efficient and short semi-group argument. In order
to make the connection, let us divide space-time into squares $Q$ of side-length
$h$ (no parabolic scaling needed here), which we think of being small. Assume that we are dealing with a function
$F$ of the space-time white noise $\xi$ that depends on $\xi$ only through the average of $\xi$ on the cubes $Q$;
which amounts to saying that $F$ only depends on $\{\xi_Q\}_{Q}$, where $\xi_Q:=\frac{1}{h}\int\xi dxdt$
(any reasonable function $F$ can be approximated by such functions $F_h$ for $h\downarrow0$).
The reason for using this normalization by the {\it square-root} of the space-time volume $h^2$
is that the application $\xi\mapsto\{\xi_Q\}_Q$ pushes
the space-time white-noise ensemble $\langle\cdot\rangle$ into the normal Gaussian
ensemble $\langle\cdot\rangle_h$. In particular $\langle F\rangle=\langle F\rangle_h$
and $\langle\exp(\frac{1}{C}F^2)\rangle=\langle\exp(\frac{1}{C}F^2)\rangle_h$.
Hence by the discrete theory, we have concentration of measure provided
we have a uniform bound on the squared Euclidean (rather Hilbertian)
norm $|\nabla_hF|^2:=\sum_{Q}(\frac{\partial F}{\partial\xi_Q})^2$
of the (infinite-dimensional) vector of partial derivatives. Therefore it remains to argue
that $|\nabla_h F|^2$ is dominated by the carr\'e-du-champs $|\nabla F|^2$
of the continuum Malliavin derivative. By definition (\ref{i.a}) of the latter we have
\begin{equation}\label{l014}
\lim_{\epsilon\downarrow 0}\frac{1}{\epsilon}(F(\xi+\epsilon\delta\xi)-F(\xi))
\le|\nabla F|\Big(\int(\delta\xi)^2dxdt\Big)^\frac{1}{2}
\end{equation}
for any field $\delta\xi$, hence in particular for a field $\delta\xi$ which is piecewise
constant on the cubes. More precisely, we may assume that $\delta\xi$ is of the
form $\delta\xi_{|Q}=\frac{1}{h}\delta\xi_Q$ for some $\{\delta\xi_Q\}_Q$
so that $\frac{1}{h}\int_Q\delta\xi dxdt=\delta\xi_Q$.
Because of this normalization, the l.\ h.\ s.\ of (\ref{l014}) turns into
$\sum_{Q}\frac{\partial F}{\partial\xi_Q}\delta\xi_Q$ by definition of the partial derivatives,
whereas the r.\ h.\ s.\ turns into $(\sum_{Q}\delta\xi_Q^2)^\frac{1}{2}$,
so that by the arbitrariness of $\{\delta\xi_Q\}_Q$,
(\ref{l014}) indeed implies $|\nabla_h F|\le|\nabla F|$ (in fact, there is equality).

\medskip

We start with the first half of the proof, that is, the bound (\ref{l04}) on the expectation.
In fact, we shall establish that
\begin{equation}\nonumber
\langle E^2(u,R)\rangle\lesssim 1,
\end{equation}
provided the scale $R\sim 1$ is sufficiently large (larger than a constant only depending on $\lambda$).
This indeed implies (\ref{l04}) since by definition of the average $\int_{-R^2}^0\int\eta_R\cdot dxdt$,
$E^2(u,1)\le R^3 E^2(u,R)$, where the power three represents the parabolic dimension.
By the scale invariance (\ref{t1}) \& (\ref{t2}), we might as well show
\begin{equation}\label{l05}
\langle E^2(u,1)\rangle\lesssim 1,
\end{equation}
provided the massive term is sufficiently strong, that is, $T\sim 1$ is sufficiently small.
In fact, it will be convenient for the upcoming calculation
to replace the exponential cut-off $\eta$ by $\tilde\eta^2$,
where $\tilde\eta$ is a smoothened version of $\eta_2$, to fix ideas
\begin{equation}\label{l08}
\tilde\eta(x):=\exp(-\frac{1}{2}\sqrt{x^2+1}).
\end{equation}
In order to establish (\ref{l05}), we will use a martingale argument based on the stochastic
(partial) differential equation with (nonlinear) damping
\begin{equation}\label{l06}
\partial_t u=-(\frac{1}{T}u+(-\partial_x^2)\pi(u))+\xi.
\end{equation}
As is constitutive for a martingale argument, we shall monitor a symmetric and
semi-definite expression, in our case $\int\tilde\eta u(1-\partial_x^2)^{-1}\tilde\eta udx$,
where we use physicist's notation in the sense that an operator, here $(1-\partial_x^2)^{-1}$,
acts on everything to its right, here the product $\tilde\eta u$. This quadratic expression,
which amounts to a version of the $\dot H^{-1}$-norm that is localized (thanks to the inclusion 
of $\tilde\eta$) and endowed with an infra-red cut-off (the effect of the $1$ in $(1-\partial_x^2)^{-1}$),
is motivated by the $\dot H^{-1}$ contraction
principle, a well-known property of the deterministic versions of (\ref{i.1}); in this language,
we monitor here the (modified) $\dot H^{-1}$ distance to the trivial solution $u=0$.
In general terms, the time derivative of such quadratic expression under a stochastic equation
comes in three contributions: the contribution solely of the deterministic r.\ h.\ s.\ of (\ref{l06}),
the contribution solely from the stochastic r.\ h.\ s.\ $\xi$, and a mixed contribution.
In this set-up, the space-time white noise $\xi$ is viewed as
a white noise in time with a spatial (and thus infinite-dimensional) covariance structure expressing
white noise in space. The mixed contribution is a martingale, and thus vanishes when taking
the expectation: This cancellation can best be understood when considering a time discretization of (\ref{l06}) that
is explicit in the drift $-(\frac{1}{T}u+(-\partial_x^2)\pi(u))$ (of course, an explicit time discretization
is not well-posed for an infinite dimensional dynamical system coming from a {\it parabolic} equation, so
one better combines it in one's mind with a spatial discretization). 
The contribution which solely comes from $\xi$ is 
the so-called quadratic variation, and its expectation
can be computed based on the operator defining the
quadratic expression, here $\tilde\eta(1-\partial_x^2)^{-1}\tilde\eta$,
and the spatial covariance structure of the noise (provided it is white in time).
Since the spatial covariance structure is the one coming from (spatial) white noise,
it is given by the integral of the diagonal of the kernel (i.\ e.\ the trace-norm of the operator).
In case of $\tilde\eta(1-\partial_x^2)^{-1}\tilde\eta$, the kernel is given
by $\tilde\eta(x)\frac{1}{2}\exp(-|x-y|)\tilde\eta(y)$. Hence the expectation of the quadratic variation
is given by $\int\frac{1}{2}\tilde\eta^2dx$. Altogether, the martingale argument thus yields 
%
\begin{eqnarray}\label{l210}
\lefteqn{\frac{d}{dt}\frac{1}{2}\Big\langle\int\tilde\eta u(1-\partial_x^2)^{-1}\tilde\eta udx
\Big\rangle}\nonumber\\
&=&-\Big\langle\int\tilde\eta u(1-\partial_x^2)^{-1}\big(\frac{1}{T}\tilde\eta u
+\tilde\eta(-\partial_x^2)\pi(u)\big)dx
\Big\rangle+\frac{1}{2}\int\frac{1}{2}\tilde\eta^2dx.
\end{eqnarray}
We rewrite this identity as
\begin{eqnarray*}
\lefteqn{\frac{d}{dt}\exp(\frac{t}{T})\Big\langle\int\tilde\eta u(1-\partial_x^2)^{-1}\tilde\eta udx
\Big\rangle}\nonumber\\
&=&\exp(\frac{t}{T})\Big(-\Big\langle\int\tilde\eta u(1-\partial_x^2)^{-1}\big(\frac{1}{T}\tilde\eta u
+2\tilde\eta(-\partial_x^2)\pi(u)\big)dx
\Big\rangle+\int\frac{1}{2}\tilde\eta^2dx\Big),
\end{eqnarray*}
and integrate over $t\in(-\infty,0)$:
\begin{eqnarray}
\Big\langle\int_{-\infty}^0\exp(\frac{t}{T})\int\tilde\eta u(1-\partial_x^2)^{-1}
\big(\frac{1}{T}\tilde\eta u+2\tilde\eta(-\partial_x^2)\pi(u)\big)dxdt\Big\rangle
\le \frac{T}{2}\int\tilde\eta^2dx.\nonumber
\end{eqnarray}
Hence in order to arrive at (\ref{l05}),
it is enough to show that for $T\ll 1$, we have the deterministic estimate
\begin{equation}\label{l07}
\int\tilde\eta u(1-\partial_x^2)^{-1}\big(\frac{1}{T}\tilde\eta u
+2\tilde\eta(-\partial_x^2)\pi(u)\big)dx
\gtrsim\int\tilde\eta^2u^2dx.
\end{equation}

\medskip

We have a closer look at the elliptic term $\tilde\eta(-\partial_x^2)\pi(u)$
in (\ref{l07}), whose contribution would be positive by the monotonicity of $\pi$ if
it weren't for the spatial cut-off and the infra-red cut off. 
Using Leibniz' rule, we rewrite it as
(in our physicist's way of omitting parentheses)
\begin{equation}\label{p4.1}
\tilde\eta(-\partial_x^2)\pi(u)
=(1-\partial_x^2)\pi(u)\tilde\eta+2\partial_x\pi(u)\partial_x\tilde\eta
-\pi(u)(\tilde\eta-\partial_x^2\tilde\eta),
\end{equation}
where we w.\ l.\ o.\ g.\ assume that $\pi(0)=0$.
Hence by the symmetry of $(1-\partial_x^2)^{-1}$ we obtain
\begin{eqnarray*}
\lefteqn{\int\tilde\eta u(1-\partial_x^2)^{-1}\tilde\eta(-\partial_x^2)\pi(u)dx}\\
&=&\int\tilde\eta^2u\pi(u)dx-2\int(\partial_x\tilde\eta)\pi(u)\partial_x(1-\partial_x^2)^{-1}\tilde\eta udx\\
&&-\int(\tilde\eta-\partial_x^2\tilde\eta)\pi(u)(1-\partial_x^2)^{-1}\tilde\eta udx.
\end{eqnarray*}
Using that the operators $\partial_x(1-\partial_x^2)^{-\frac{1}{2}}$ and $(1-\partial_x^2)^{-\frac{1}{2}}$
have operator norm 1 w.\ r.\ t.\ to $L^2$, we deduce the inequality
\begin{eqnarray}\nonumber
\lefteqn{\int\tilde\eta u(1-\partial_x^2)^{-1}\tilde\eta(-\partial_x^2)\pi(u)dx
\;\ge\;\int\tilde\eta^2u\pi(u)dx}\nonumber\\
&&-\Big(2\big(\int(\partial_x\tilde\eta)^2\pi^2(u)dx\big)^\frac{1}{2}
+\big(\int(\tilde\eta-\partial_x^2\tilde\eta)^2\pi^2(u)dx\big)^\frac{1}{2}\Big)\nonumber\\
&&\times \left(\int\tilde\eta u(1-\partial_x^2)^{-1}\tilde\eta udx\right)^\frac{1}{2}.
\end{eqnarray}
By the monotonicity properties (\ref{i.2}) of $\pi$ and our gratuitous assumption $\pi(0)=0$, this yields
\begin{eqnarray}\nonumber
\lefteqn{\int\tilde\eta u(1-\partial_x^2)^{-1}\tilde\eta(-\partial_x^2)\pi(u)dx}\nonumber\\
&\ge&\lambda\int\tilde\eta^2u^2dx-\Big(2\big(\int(\partial_x\tilde\eta)^2u^2dx\big)^\frac{1}{2}
+\big(\int(\tilde\eta-\partial_x^2\tilde\eta)^2u^2dx\big)^\frac{1}{2}\Big)\nonumber\\
&&\times \left(\int\tilde\eta u(1-\partial_x^2)^{-1}\tilde\eta udx\right)^\frac{1}{2}.\nonumber
\end{eqnarray}
Our smoothing out of the exponential cut-off function, cf.\ (\ref{l08}), has the sole purpose
of making sure that
\begin{equation}\label{p4.2}
|\partial_x\tilde\eta|+|\partial_x^2\tilde\eta|\lesssim\tilde\eta,
\end{equation}
so that we obtain by Young's inequality for the elliptic term,
\begin{eqnarray}\nonumber
\int\tilde\eta u(1-\partial_x^2)^{-1}\tilde\eta(-\partial_x^2)\pi(u)dx
&\ge&\frac{1}{C}\int\tilde\eta^2u^2dx
-C\int\tilde\eta u(1-\partial_x^2)^{-1}\tilde\eta udx.\nonumber
\end{eqnarray}
We thus see that thanks to the massive term, (\ref{l07}) holds for $T\ll 1$.

\medskip

We now turn to the second half of the proof, the estimate of the carr\'e-du-champs (\ref{l01}).
We first argue that (\ref{l01}) follows from the deterministic estimate
\begin{equation}\label{l02}
E^2(\delta u,1)\lesssim \int(\delta\xi)^2dxdt,
\end{equation}
where $\delta u$ and $\delta\xi$ are related via
\begin{equation}\label{l016}
\delta u+\partial_t\delta u-\partial_x^2(a\delta u)=\delta\xi
\end{equation}
with $a=\pi'(u)$. Indeed, we note that by duality w.\ r.\ t.\ to
the inner product $(g,f)\mapsto \int_{-1}^0\int\eta gfdxdt$,
\begin{eqnarray}\label{l015}
\lefteqn{E(u,1)=\sup\Big\{\;E(u,f):=\int_{-1}^0\int\eta u fdxdt\;\Big|}\nonumber\\
&&\;\int_{-1}^0\int\eta f^2dxdt=1,\;\;{\rm supp}f\subset(-1,0)\times\mathbb{R}\;\Big\}.
\end{eqnarray}
By the chain rule for the Malliavin derivative we thus obtain
\begin{equation}\nonumber
|\nabla E(\cdot,1)|\le\sup_f|\nabla E(\cdot,f)|,
\end{equation}
where the supremum runs over the set implicitly defined in (\ref{l015}), so that it
is enough to show for a fixed $f$
\begin{equation}\nonumber
|\nabla E(u,f)|^2\lesssim 1.
\end{equation}
By definition of the carr\'e-du-champs of the Malliavin derivative in case of
the {\it linear} functional $u\mapsto E(u,f)$, cf.\ (\ref{i.a}), this amounts to showing
\begin{equation}\nonumber
\int_{-1}^0\int\eta \delta u fdxdt
\lesssim 1,
\end{equation}
where the infinitesimal perturbation $\delta u$ of the solution is related
to the infinitesimal perturbation $\delta\xi$ of the noise via (\ref{l016}).
By the characterizing properties of the $f$'s, cf.\ (\ref{l015}), this estimate
in turn amounts to establishing (\ref{l02}).

\medskip

We now turn to the proof of the deterministic estimate (\ref{l02}). To ease notation and
make the connection to Proposition \ref{P4}, we rephrase (and strengthen) the goal: For
$w$ and $f$ related via
\begin{equation}\label{l010}
w+\partial_tw-\partial_x^2(aw)=f,
\end{equation}
with uniformly elliptic coefficient field $a$ in the sense of (\ref{p21}), we seek the estimate
\begin{equation}\label{l09}
\int_{-\infty}^0\int w^2dxdt\lesssim\int_{-\infty}^0\int f^2dxdt.
\end{equation}
Like for (\ref{l05}),
our Ansatz for (\ref{l09}) is motivated by the $\dot H^{-1}$-contraction principle. Again, we consider a
version of $\dot H^{-1}$-norm with ultra-red cut-off, but this time without cut-off function $\eta$,
namely $\int w(L^{-2}-\partial_x^2)^{-1}wdx$, where the length scale $L$ for the ultra-red cut-off
will be chosen later. We obtain from the equation (\ref{l09})
\begin{eqnarray*}
\lefteqn{\frac{d}{dt}\frac{1}{2}\int w(L^{-2}-\partial_x^2)^{-1}wdx}\nonumber\\
&=&
-\int w(L^{-2}-\partial_x^2)^{-1}(w-f+(-\partial_x^2)(aw))dx\\
&=&
-\int w(L^{-2}-\partial_x^2)^{-1}(w-f-L^{-2}aw)dx-\int aw^2dx.
\end{eqnarray*}
We apply Cauchy-Schwarz' inequality and use the uniform ellipticity of $a$,
cf.\ (\ref{p21}), to obtain the estimate
\begin{eqnarray*}
\lefteqn{\frac{d}{dt}\frac{1}{2}\int w(L^{-2}-\partial_x^2)^{-1}wdx}\nonumber\\
&\le&
-\int w(L^{-2}-\partial_x^2)^{-1}wdx-\lambda\int w^2dx\\
&&+\Big(\int ((L^{-2}-\partial_x^2)^{-1}w)^2dx\Big)^\frac{1}{2}
 \Big(\big(\int f^2dx\big)^\frac{1}{2}+L^{-2}\big(\int w^2dx\big)^\frac{1}{2}\Big).
\end{eqnarray*}
Thanks to the operator inequality $(L^{-2}-\partial_x^2)^{-1}\le L(L^{-2}-\partial_x^2)^{-\frac{1}{2}}$
we have
\begin{equation}\nonumber
\Big(\int ((L^{-2}-\partial_x^2)^{-1}w)^2dx\Big)^\frac{1}{2}
\le L\big(\int w(L^{-2}-\partial_x^2)^{-1}wdx\Big)^\frac{1}{2},
\end{equation}
so that we may absorb the term
$(\int ((L^{-2}-\partial_x^2)^{-1}w)^2dx)^\frac{1}{2}L^{-2}(\int w^2dx)^\frac{1}{2}$ 
by Young's inequality for $L\gg 1$, obtaining
\begin{eqnarray*}
\frac{d}{dt}\int w(L^{-2}-\partial_x^2)^{-1}wdx
&\le&
-\frac{1}{C}\int w^2dx+CL^2\int f^2dx.
\end{eqnarray*}
Integration in time yields (\ref{l09}).

\bigskip

{\sc Proof of Proposition \ref{P4}}.
 We first note that (\ref{p24}) follow easily from (\ref{p43}): 
by translation invariance
(\ref{p43}) also holds with $\eta$ replaced by the shift $\eta^y$, summation over
$y\in\mathbb{Z}$ gives (\ref{p24}). We next note that w.\ l.\ o.\ g.\ we
may assume $h=0$, since we may rewrite (\ref{i.5}) as
$\partial_t(w-h)-\partial_x^2(a(w-h))=\partial_x^2(g+ah)$. In this form \eqref{p43} follows 
from \eqref{p43A} which we will proceed to show now.
The proof of this proposition
is very close to the deterministic part of the proof of Lemma \ref{L0};
in fact, it might be seen as an infinitesimal version of it. Like there,
we substitute $\eta$ by $\tilde\eta^2$, cf.\ (\ref{l08}), and start from monitoring
the localized $H^{-1}$-norm of $w$ with infra-red cut-off:
\begin{equation}\nonumber
\frac{d}{dt}\frac{1}{2}\int\tilde\eta w(1-\partial_x^2)^{-1}\tilde\eta w dx
=-\int\tilde\eta w(1-\partial_x^2)^{-1}\tilde\eta(-\partial_x^2)(aw+g)dx.
\end{equation}
As in (\ref{p4.1}), we write
\begin{equation}\nonumber
\tilde\eta(-\partial_x^2)(aw+g)
=(1-\partial_x^2)(aw+g)\tilde\eta+2\partial_x(aw+g)\partial_x\tilde\eta
-(aw+g)(\tilde\eta-\partial_x^2\tilde\eta),
\end{equation}
which yields
\begin{eqnarray*}\nonumber
\frac{d}{dt}\frac{1}{2}\int\tilde\eta w(1-\partial_x^2)^{-1}\tilde\eta w dx
&=&-\int\tilde\eta w\tilde\eta(aw+g)dx\\
&&-2\int\tilde\eta w(1-\partial_x^2)^{-1}\partial_x(aw+g)\partial_x\tilde\eta dx\\
&&+\int\tilde\eta w(1-\partial_x^2)^{-1}(aw+g)(\tilde\eta-\partial_x^2\tilde\eta) dx.
\end{eqnarray*}
Using symmetry and boundedness properties of $(1-\partial_x^2)^{-1}$, 
and the estimates (\ref{p4.2}) on our mollified exponential cut-off $\tilde\eta$, 
the two last terms are estimated as
\begin{eqnarray*}
\lefteqn{\int\tilde\eta w(1-\partial_x^2)^{-1}\partial_x(aw+g)\partial_x\tilde\eta dx}\\
&\lesssim&\Big(\int\tilde\eta w(1-\partial_x^2)^{-1}\tilde\eta wdx
\int\tilde\eta^2(aw+g)^2dx\Big)^\frac{1}{2}
\end{eqnarray*}
and
\begin{eqnarray*}
\lefteqn{-\int\tilde\eta w(1-\partial_x^2)^{-1}(aw+g)(\tilde\eta-\partial_x^2\tilde\eta) dx}\\
&\lesssim&\Big(\int\tilde\eta w(1-\partial_x^2)^{-1}\tilde\eta wdx
\int\tilde\eta^2(aw+g)^2dx\Big)^\frac{1}{2}.
\end{eqnarray*}
Hence we obtain by the uniform ellipticity (\ref{p21}) of $a$ together
with the triangle inequality to break up $aw+g$ and Young's inequality
\begin{eqnarray*}
\lefteqn{\frac{d}{dt}\int\tilde\eta w(1-\partial_x^2)^{-1}\tilde\eta wdx
+\frac{1}{C}\int\tilde\eta^2 w^2dx}\\
&\le&C\Big(\int\tilde\eta w(1-\partial_x^2)^{-1}\tilde\eta w dx+\int\tilde\eta^2 g^2dx\Big),
\end{eqnarray*}
which we rewrite as
\begin{eqnarray*}
\lefteqn{\frac{d}{dt}\exp(-Ct)\int\tilde\eta w(1-\partial_x^2)^{-1}\tilde\eta wdx}\\
&&+\frac{1}{C}\exp(-Ct)\int\tilde\eta^2 w^2dx
\;\le\;C\exp(-Ct)\int\tilde\eta^2 g^2dx.
\end{eqnarray*}
The desired estimate \eqref{p43A} follows.
%

\bigskip

{\sc Proof of Proposition \ref{P5}.} 
We will again work with a smooth version $\te$ of the cutoff $\sqrt{\eta}$, setting  
\begin{equation*}
\te(x) = \frac{1}{2} \exp \big( - \frac12 \sqrt{|x|^2+1}   \big)  \qquad \te_r(x) = \te (x/r).
\end{equation*}
Note in particular that this time the  cut-off $\te_r$ at scale $r$ is not normalised to preserve the $L^1$ norm.
We will establish \eqref{i.10} in the form
\begin{align}\label{P4-0}
\int_{-r^2}^0 \int \te_r^2 v^2 dx dt \lesssim r^{1 + 2 \an}\int (\te v) (1- \px^2)^{-1} (\te v) dx \big|_{t=-1}.
\end{align}
for some $\an>0$.

\medskip
We start by defining the auxiliary function $V = -\partial_x (1-\partial_x^2 )^{-1}(\te v)$, which we think 
of as a localised version of the anti-derivative of $v$. We claim that
\begin{align}\label{P4-1}
\partial_t V - \partial_x (a \px V) = \partial_x g +f \,
\end{align}
where
\begin{align*}
g &:  =  a  (1- \px^2)^{-1} (\te v),  \\
f &:=2 \px^2  (1-\partial_x^2 )^{-1} (  a (\px \te)   v )   - \px (1-\partial_x^2 )^{-1} (a (\te + \partial_x^2 \te) v) .
\end{align*}
To see \eqref{P4-1} we first observe that
\begin{align}
\label{P4-2}
\partial_t V &= -\partial_x (1-\partial_x^2 )^{-1}(\te  \partial_t v) =  -\partial_x (1-\partial_x^2 )^{-1}(\te  \partial_x^2 (a v)). 
\end{align}
Then we write
\begin{align}
\notag
\te  \partial_x^2 (a v) &=  \partial_x^2 (\te a v)  - 2 (\px \te) \partial_x (a v) - (\partial_x^2 \te ) av \\
\notag
&=  \partial_x^2 (\te a v)  - 2 \px \big(  (\px \te)  a v \big) +  (\partial_x^2 \te ) av.
\end{align}
Plugging this into \eqref{P4-2} we get
\begin{align}
\notag
\partial_t V  &= -\partial_x (1-\partial_x^2 )^{-1}\Big[  \partial_x^2 (a \te  v)  - 2 \px \big( a  (\px \te)   v \big) +  a (\partial_x^2 \te ) v \Big]  \\
\notag
& =  \px  (a \te  v) - \px (1-\partial_x^2 )^{-1}\Big[  a \te  v - 2 \px \big(a   (\px \te)   v \big) + a (\partial_x^2 \te ) v \Big] \\
\label{P4-2A}
& =  \px  (a \te  v) + f .
\end{align}
On the other hand we have
\begin{align}\label{P4-2B}
\px V= -\px^2 (1- \px^2)^{-1} (\te v) = \te v - (1- \px^2)^{-1} (\te v) ,
\end{align}
which together with \eqref{P4-2A} implies that 
\begin{align*}
\partial_t V  & =  \px  ( a \px V) +\px \big(a  (1- \px^2)^{-1} (\te v) \big) + f .
\end{align*}
So \eqref{P4-1} follows.

\medskip

Our next step is to derive a suitable version of Caccioppoli's estimate 
 on scale $r \ll 1$ for \eqref{P4-1}. We can write for any $c \in \mathbb{R}$
\begin{align}
\notag
&\partial_t  \int \te_r^2 \frac12 (V-c)^2 \, dx \\
\notag
&\qquad = -\int \px \big( \te_r^2 (V-c) \big) \Big[ a \px V  + g \Big] \; dx + \int  \te_r^2 (V-c) f \; dx  \\
\label{P4-3}
&\qquad = -\int \Big( 2 ( \px \te_r) \te_r (V-c)  +  \te_r^2 \px V \big) \Big[ a \px V  + g \Big] \; dx  + \int  \te_r^2 (V-c) f \; dx  .
\end{align}
We treat the terms on the right hand side of this expression one by one. First we get
\begin{align}
\notag
&- \int \Big( 2 ( \px \te_r) \te_r (V-c)  +  \te_r^2 \px V \big)  a \px V \, dx \\
\notag
&\quad  \leq- \lambda \ \int \te_r^2 (\px V)^2 \, dx +2 \Big(\int (\px \te_r)^2 (V-c)^2 dx \Big)^{\frac12}\Big(\int \te^2_r (\px V)^2 \; dx \Big)^{\frac12} \\
\label{P4-3A}
&\quad  \leq- \frac{1}{C} \int \te_r^2 (\px V)^2 \, dx +  \frac{C}{r^2} \int \te_r^2 (V-c)^2 dx  ,
\end{align}
where in the last line we have absorbed the second term involving $\px V$ in the first one and used the point-wise estimate $|\px \te_r| \lesssim \frac{1}{r} \te_r$.
For the second term on the right hand side of \eqref{P4-3} we get
\begin{align}
\label{P4-4}
&-\int \Big( 2 ( \px \te_r) \te_r (V-c)  +  \te_r^2 \px V \big)  g  \; dx \\
\notag
 &  \lesssim \Big( \frac{1}{r^2}\int  \te^2 (V-c)^2 dx \Big)^{\frac12}  \Big( \int \te_r^2 g^2 \, dx\Big)^{\frac12} + \Big( \int \te_r^2 (\px V)^2 dx \Big)^{\frac12}  \Big( \int \te_r^2 g^2 \, dx\Big)^{\frac12}  .
\end{align}
The integral involving $(V-c)$ and the integral involving $\px V$ can be absorbed into the terms on the right hand side of \eqref{P4-3A}. The only contribution that is left  
is the integral $\int \te_r^2 g^2 \, dx$.
Finally, for the last term on the right hand side of \eqref{P4-2} we get
\begin{align}
\notag
 \int  \te_r^2 (V-c) f \; dx  \lesssim \Big( \sup_x \te_r |V-c| \Big) \int \te_r |f| \, dx.   
\end{align}
Summarising, we obtain 
\begin{align}
\notag
&  \frac{d}{dt}\int  \te_r^2 \frac12 (V-c)^2 \, dx +\frac{1}{C} \int \te_r^2 (\px V)^2 \, dx dt \\
\label{P4-5}
&\quad \lesssim \frac{1}{r^2} \int \te_r^2 (V-c)^2 dx +    \int \te_r^2 g^2 \, dx + \Big( \sup_x \te_r |V-c| \Big) \int \te_r |f| \, dx .
\end{align}
As a next step we integrate this estimate in time. To this end, let  $\zeta\colon \mathbb{R} \to \mathbb{R}$ be non-negative, non-decreasing  such that $\zeta = 0$ on $(-\infty,-2)$, $\zeta=1$ on $(-1,\infty)$ and with $ \zeta' \leq 2$ on $\mathbb{R}$.
We set $\zeta_r(t) = \zeta(t/r^2)$. Then integrating \eqref{P4-5} against $\zeta_r$  we get
\begin{align}
\notag
&\int_{-\infty}^0   \frac{d}{dt}  \int \zeta_r  \te_r^2 \frac12 (V-c)^2 \, dxdt + \int_{-r^2}^0 \int \te_r^2 (\px V)^2 \, dx dt\\
\notag
&\quad \lesssim \frac{1}{r^2} \int_{-2r^2}^0 \int \te_r^2 (V-c)^2 dxdt +   \int_{-2r^2 }^0 \int \te_r^2 g^2 \, dxdt \\
\label{P4-5A}
& \qquad + \Big( \sup_{(-2r^2,0) \times \mathbb{R}}  \te_r |V-c| \Big) \int_{-2r^2}^0\int \te_r |f| \, dxdt ,
\end{align}
where we have absorbed the term $\int_{-1}^0 \int  (\partial_t  \zeta_r)    \te_r^2 \frac12 (V-c)^2 \, dx dt$ in the first term on the right hand side.
We proceed by bounding the first term on the right hand side 
\begin{align*}
 \frac{1}{r^2}\int_{-2r^2}^0 \int \te_r^2 (V-c)^2 dxdt &\leq  \Big( \sup_{(-2r^2,0) \times \mathbb{R}}  \te_r (V-c)^2 \Big) \frac{1}{r^2}\int_{-2r^2}^0 \int \te_r dx dt  \\
 &\lesssim  r \sup_{(-2r^2,0) \times \mathbb{R}} \te_r (V-c)^2 , 
\end{align*}
the second term by 
\begin{align*}
 \int_{-2r^2}^0  \int \te_r^2 g^2 \, dxdt \lesssim r^3  \sup_{(-1,0) \times \mathbb{R}}   g^2
\end{align*}

and the last term by
\begin{align*}
 &\Big( \sup_{(-2r^2,0) \times \mathbb{R}}  \te_r |V-c| \Big) \int_{-2r^2}^0\int \te_r |f| \, dxdt\\
 &\qquad  \lesssim r \Big( \sup_{(-2r^2,0) \times \mathbb{R}}  \te_r |V-c| \Big)^2 + \frac{1}{r} \Big( \int_{-2r^2}^0\int \te_r |f|  \; dx\Big)^2\\
&\qquad   \lesssim r \Big( \sup_{(-2r^2,0) \times \mathbb{R}}  \te_r |V-c| \Big)^2 +   r^2 \int_{-1}^0\int  |f|^2  \; dx.
\end{align*}

\medskip
Inserting these estimates in \eqref{P4-5A} we arrive at 
\begin{align}
\notag
&\int_{-r^2}^0 \int \te_r^2 (\px V)^2 \, dxdt\\
& \lesssim 
r \sup_{(-2r^2,0) \times \mathbb{R}} \te_r (V-c)^2 + r^3  \sup_{(-1,0) \times \mathbb{R}} g^2 + r^2 \int_{-1}^0 \int f^2  \; dx dt ,\label{P4-6}
\end{align}
where we have dropped the non-negative term $\int \te_r^2 \frac12 (V-c)^2 \, dx|_{t=0} $ on the left hand side.
\medskip

In the following crucial step we will use the De Giorgi-Nash Theorem to obtain a slightly larger power of  $r$ in the first term on the right hand side of \eqref{P4-6}. More precisely, 
we will use the estimate 
\begin{align*}
[V]_{\an, (-\frac12,0) \times (-\frac12, \frac12)} \lesssim& \Big( \int_{-1}^0 \int_{-1}^1 V^2  dx dt \Big)^{\frac12} +  \Big( \int_{-1}^0 \int_{-1}^1 f^2  dx dt \Big)^{\frac12} \\
&+ \sup_{(-1,0)\times (-1,1)}|g| ,
\end{align*}
where $[V]_{\an}$ denotes the parabolic $\an$ H\"older norm of $V$ (defined as in \eqref{DefHol}). 
We refer to Theorem 8.1 and Theorem 10.1 
in \cite{LSU}; Theorem 8.1
gives local control of $\sup_{t,x}|V|$ in terms of $(\int V^2dxdt)^\frac{1}{2}$, and Theorem 10.1
gives local control of $[V]_{\an}$ in terms of $\sup_{t,x}|V|$. Our control of the right hand side
$g$ and $f$ in (\ref{P4-1}) through the local norms $\sup_{t,x}|g|$ and $(\int f^2dxdt)^\frac{1}{2}$ 
is well within the allowed range, cf.\ (7.1) \& (7.2) in \cite{LSU}, which in one
space dimension can deal with control of $\int\int (|g|^2+f)^2dxdt$ ($n=1$, $q=r=2$, and $\kappa=\frac{1}{4}$
in the notation of this reference).

\medskip

We assemble the last estimate to 
\begin{align*}
\sup_{k \in \mathbb{Z} }[V]_{\an, (-\frac12,0) \times (k -\frac12, k+\frac12)} \lesssim& \Big( \int_{-1}^0 \int_{-\infty}^\infty V^2  dx dt \Big)^{\frac12} +  \Big( \int_{-1}^0 \int_{-\infty}^{\infty} f^2  dx dt \Big)^{\frac12} \\
&+ \sup_{(-1,0)\times (-\infty,\infty)}|g| .
\end{align*}
Going back to \eqref{P4-6} we choose $c = V(0,0)$ and write
\begin{align*}
&\sup_{(t,x)  \in (-2r^2,0) \times \mathbb{R}} \te_r (V(t,x)-V(0,0) )^2\\
& \lesssim  
 \sum_{k \in \mathbb{Z}}  e^{-|k| } \sup_{(t,x)  \in (-2r^2,0) \times [r(k - \frac{1}{2}), r(k+ \frac{r}{2})] } (V(t,x)-V(0,0) )^2 \\
& \lesssim r^{2\an} \;  \sup_{k \in \mathbb{Z} }[V]_{\an, (-\frac12,0) \times (k -\frac12, k+\frac12)}.
\end{align*}
Combining these estimates we obtain the following inverse H\"older inequality
\begin{align}
\notag
&\int_{-r^2}^0 \int \te_r^2 (\px V)^2 dx dt \\
\label{P4-7}
& \lesssim r^{1+ 2 \an} \Big[ \int_{-1}^0 \int  V^2 + \int_{-1}^0 \int f^2 \, dx dt \, + \sup_{(-1,0) \times \mathbb{R}} g^2 \Big]
\end{align}
valid for $r \ll 1$. 

\medskip

It remains to control the functions $f$ and $g$. 
For $f$ we have by definition
\begin{align*}
\int_{-1}^0 \int f^2 dx dt  &\lesssim \int_{-1}^0 \int \Big(  \px^2  (1-\partial_x^2 )^{-1} (  a (\px \te)   v ) \Big)^2  dx dt\\
& \qquad  +  \int_{-1}^0 \int  \Big( \px (1-\partial_x^2 )^{-1}  (\te + \partial_x^2 \te)(a v) \Big)^2 dxdt.
\end{align*}
Using the boundedness of $ \px^2  (1-\partial_x^2 )^{-1}$ on $L^2$, as well as $|a |\leq1$ and $ |\px \te| \lesssim \te$ we can see that the first integral is bounded by $\int (   \te  v )^2 dx$. Using the boundedness of $  \px  (1-\partial_x^2 )^{-1}$ on $L^2$ and the point-wise bound $|\te + \partial_x^2 \te|\lesssim \te$ we bound the second integral by the same quantity.  
For $g$ we write for any $t$ using the embedding $H^1 \hookrightarrow L^\infty$
\begin{align*}
\sup_{x \in \mathbb{R}} g^2
& \lesssim \sup_x \big(  \;(1- \px^2)^{-1} (\te v) \big)^2\\
& \lesssim \int \big( \px (1- \px^2)^{-1} (\te v) \big)^2 \, dx + \int \big(  (1- \px^2)^{-1} (\te v) \big)^2 \, dx  \\
& =   \int  (1- \px^2)^{-1} (\te v) \;\big( -\px^2 (1- \px^2)^{-1} (\te v)\big)  \, dx  +  \int \big(  (1- \px^2)^{-1} (\te v) \big)^2 \, dx \\
&\lesssim \int (\te v)    (1- \px^2)^{-1} (\te v)    \, dx  ,
\end{align*}
where in the last step we have used the boundedness of $(1 - \px^2)^{-\frac12}$ on $L^2$.
\medskip

\medskip

Summarising these bounds and plugging in the identity  \eqref{P4-2B} which expresses $\px V$ as $\te v$ and a higher order term we obtain
\begin{align}
\notag
&\int_{-r^2}^0 \int \te_r^2 (\te v)^2 dx dt \\
\notag
& \lesssim r^{1+ 2 \an} \Big[ \int_{-1}^0 \int (\px (1-\px^2)^{-1} \te v)^2 dx dt + \int_{-1}^0 \int (\te v)^2 dx dt  \\
&\qquad+ \sup_{t \in (-1,0)} \int (\te v)    (1- \px^2)^{-1} (\te v)    \, dx \Big] 
 + \int_{-r^2}^0 \int \te_r^2 \big((1 - \partial_x^2)^{-1} (\te v)\big)^2 dx dt .
\label{P4-7a}
\end{align}
The first term on the right hand side is bounded by the second term due to the boundedness of $\px(1 - \px^2)^{-1}$ on $L^2$.
For the last term on the right hand side we can write, using once more the embedding $H^{1} \hookrightarrow L^\infty$
\begin{align*}
\int_{-r^2}^0 \int \te_r^2 \big((1 - \partial_x^2)^{-1} (\te v)\big)^2 dx dt \lesssim r^3  \sup_{t\in(-1,0)} \int (\te v)    (1- \px^2)^{-1} (\te v)    \, dx,
\end{align*}
so this term can be absorbed in the second term on the right hand side of \eqref{P4-7a}. 
In this form  the estimate trivially holds $r$ which are bounded away from $0$, so that we can conclude that  it holds for  $r \leq 1$. 
Then the desired estimate \eqref{P4-0} follows from
estimate \eqref{p43A} in Proposition~\ref{P4}. We finally note that  we may replace the cut-off $\te$ by $\sqrt{\eta}$ 
because the kernel of $(1-\px^2)^{-1}$ is non-negative.

\bigskip

\bigskip
{\sc Proof of Proposition \ref{P1}}.
For conciseness, we ignore the massive term in (\ref{i.1}).
The main object of this proposition is $\delta u:=u^h-u$, where $u^h(t,x)=u(t,x+h)$ denotes
a spatial shift of the stationary solution of (\ref{i.1}). We note that $\delta u$
satisfies the formally linear equation
\begin{equation}\label{p110}
\partial_t\delta u-\partial_x^2(a_h\delta u)=(\partial_t-a_0\partial_x^2)\delta g,
\end{equation}
where we introduced the coefficient field
\begin{equation}\label{p111}
a_h=\int_0^1\pi'(\sigma u^h+(1-\sigma)u)d\sigma,
\end{equation}
which by (\ref{i.2}) is uniformly elliptic in the sense of (\ref{p21}),
and we have set $\delta g:=g^h-g$, where $g$ is defined via the linear version of (\ref{i.1})
\begin{equation}\label{p120}
\partial_tg-a_0\partial_x^2g=\xi\quad\mbox{for}\;\;t\in(-1,0),\quad
g=0\;\;\mbox{for}\;t=-1,
\end{equation}
cf.\ Lemma \ref{L3}, with a constant coefficient $a_0\in[\lambda,1]$ to be chosen below.

\medskip

We start with the main deterministic ingredient for Proposition \ref{P1},
which we need to go from scales of order one to scales of order $r\ll 1$
in an $L^2$-averaged H\"older modulus of continuity. It is given by the estimate
\begin{align}
\notag
&\av_{-r^2}^0\int\eta_r(\delta u)^2dxdt 
 \lesssim r^{ 2\an} \big( \frac{h}{r} \big)^2e^{|h|} \Dprime \\
\label{p12}
&  \qquad \qquad \qquad + \Big(1+\frac{1}{r^3}\int_{-1}^0\int\eta(a_h-a_0)^2dxdt\Big)\sup_{(t,x)\in(-1,0)\times\mathbb{R}}\eta(\delta g)^2, 
\end{align}
for $\an>0$ from Proposition~\ref{P5}, which we shall establish for all $r\ll 1$.

\medskip

To this purpose, we split the 
solution $\delta u=\delta g+v+w$, where $v$ is defined through the 
initial value problem with homogeneous right hand side
\begin{equation}\label{p121}
\partial_tv-\partial_x^2 (a_h v)=0\;\;\mbox{for}\;t\in(-1,0),\quad
v=\delta u\;\;\mbox{for}\;t=-1,
\end{equation}
and where $w$ is defined through the initial value problem with homogenous initial data
\begin{equation}\label{p18}
\partial_tw-\partial_x^2(a_hw)=\partial_x^2((a_h-a_0)\delta g)\;\;\mbox{for}\;t\in(-1,0),\quad
w=0\;\;\mbox{for}\;t=-1.
\end{equation}
Taking the sum of (\ref{p120}), \eqref{p121} and (\ref{p18}), and comparing with (\ref{p110}), we see
that this indeed gives $\delta u=\delta g+v+w$.

\medskip

We first address $v$. From the estimate (\ref{i.10}) in Proposition \ref{P5}, we learn that there exists an 
$\an>0$ (depending only on $\lambda$) such that
\begin{equation}\label{p15a}
\av_{-r^2}^0\int \eta_r v^2 dx dt\lesssim
 r^{-2+ 2\an} \int \eta_{2}\delta u (1-\px^2)^{-1} \eta_{2}\delta u dx_{|t=-1}.
\end{equation}
 By Leibniz' rule an its discrete form of 
 \begin{align*}
 \eta_2 \delta u& = \eta_ 2 \delta (u-c) = \delta (\eta_2(u-c)) - \delta \eta_2 (u^h-c)\\
 & = \delta (\eta_2(u-c)) - \big(\delta \eta_2^{-h}(u-c) \big)^h
 \end{align*}
and the triangle inequality we can write the integral on the right hand side of \eqref{p15a} as 
\begin{align*}
\int &\eta_{2}\delta u (1-\px^2)^{-1} \eta_{2} \delta u dx \\
&\lesssim \int (\delta V )^2  dx+  \int \big( (1-\px^2)^{-\frac12} \delta  \eta_{2}^{-h} (u-c) \big)^2 dx ,
\end{align*}
where we have set $V := (1-\px^2)^{-\frac12} \eta_{2} (u-c)$, 
and $c$ can be chosen arbitrarily.
We use the point-wise bound $|\delta \eta_{2}^{-h}| \leq |h| e^{\frac{|h|}{2}} \eta_{2}$, the point-wise bound
$\eta_{2}^2 \lesssim \eta$
as well as the boundedness
of $(1-\px^2)^{-\frac12} $ on $L^2$ to bound the second term as follows
\begin{align*}
 \int \big( (1-\px^2)^{-\frac12} \delta  \eta_2^{-h} (u-c) \big)^2 dx \lesssim h^2 e^{|h|}  \int  \eta (u-c)^2 dx .
\end{align*}
For the first term we get
\begin{align*}
 \int (\delta V) ^2  dx &\lesssim h^2 \int (\partial_x V)^2 dx = h^2\int \big( \px  (1-\px^2)^{-\frac12} \eta_{2} (u-c) \big)^2 dx \\
 &\lesssim h^2  \int  \eta (u-c)^2 dx  ,
\end{align*}
by the boundedness on $L^2$ of $\px (1-\px^2)^{-1}$. Summarising these bounds and chosing $c = \int \eta u dx |_{t=1}$ we get
\begin{equation*}
\av_{-r^2}^0\int \eta_r v^2 dx dt\lesssim
  r^{2\an} \big(\frac{h}{r}\big)^2 e^{|h|} \Dprime  .
\end{equation*}

\medskip
We now turn to $w$.
Applying the first part (\ref{p43}) of Proposition \ref{P4} (with $\eta$ replaced by $\eta_\frac{1}{2}$), 
to \eqref{p18} we gather that
\begin{equation}\nonumber
\int_{-1}^0\int\eta_{\frac{1}{2}}w^2dxdt\lesssim\int_{-1}^0\int\eta_\frac{1}{2}(a_h-a_0)^2(\delta g)^2 dxdt,
\end{equation}
which implies for $r\le\frac{1}{2}$ (by the obvious inequality $\av_{-r^2}^0\int\eta_r\cdot dxdt
\le (\frac{R}{r})^3  \times \av_{-R^2}^0\int\eta_R\cdot dxdt$ for $r\le R$ and since $\eta_\frac{1}{2}\lesssim\eta^2$)
\begin{eqnarray}
r^3\av_{-r^2}^0\int\eta_rw^2dxdt
\lesssim
\int_{-1}^0   \int\eta(a_h-a_0)^2dxdt
\sup_{(t,x)\in(-1,0)\times\mathbb{R}} \eta (\delta g)^2.\nonumber
\end{eqnarray}
Finally, because of $\av_{-r^2}^0\int\frac{\eta_r}{\eta} dxdt\lesssim 1$ for $r\ll 1$,
we have for the last contribution $\delta g$
\begin{equation}\nonumber
\av_{-r^2}^0\int\eta_r(\delta g)^2dxdt\lesssim\sup_{(t,x)\in(-1,0)\times\mathbb{R}}\eta(\delta g)^2.
\end{equation}
Combining the four last estimates yields (\ref{p12}) for $\delta u=\delta g+v+w$.

\medskip

We now post-process (\ref{p12}) and to that purpose make the choice of
$a_0=\pi'(c)$ with $c:=\int_{-1}^0\int\eta udxdt$, so that in view of the definition (\ref{p111}) of $a_h$
and the Lipschitz continuity (\ref{i.3}) of $\pi'$
\begin{equation}\label{e:ah-vs-anull}
|a_h-a_0|\le|a_h-\pi'(u)|+|\pi'(u)-\pi'(c)|\le L(|\delta u|+|u-c|).
\end{equation}
Therefore, (after replacing $r$ by $2r$  in order to
make $\eta_{2r}$ appear, which is no problem thanks to $r\ll 1$)  (\ref{p12})  turns into
\begin{align}
\label{p114a}
&\av_{-r^2}^0\int\eta_{2r}(u^h-u)^2dxdt\\
\notag
& \qquad \lesssim r^{  2\an} \big(\frac{h}{r}\big)^2e^{|h|} \Dprime \\ 
\notag
&\qquad +\Big(1+\frac{L^2}{r^3}\int_{-1}^0\int\eta(u^h-u)^2dxdt+\frac{L^2}{r^3}D^2(u,1)\Big)
\sup_{(t,x)\in(-1,0)\times\mathbb{R}}\eta(g^h-g)^2.
\end{align}
Now we integrate in $h$ according to $\int\eta_{2r}(h)\cdot dh$. As we shall argue below, we have
for the l.\ h.\ s.\ of (\ref{p114a})
\begin{equation}\label{p115a}
\int\eta_{2r}(h)\av_{-r^2}^0\int\eta_{2r}(u^h-u)^2dxdtdh
\gtrsim\av_{-r^2}^0\int\eta_r(u-\int\eta_ru)^2dxdt.
\end{equation}
For the first term on the right hand side of \eqref{p114a} we observe that for $r \leq \frac14$
\begin{align*}
\int \eta_{2r}(h)  r^{  2\an} \big(\frac{h}{r}\big)^2e^{|h|} \Dprime d h \lesssim r^{2\an}  \Dprime.
\end{align*}
The second term on the r.\ h.\ s.\ of (\ref{p114a}) comes in form of a product of
two $h$-dependent functions we momentarily call $f_1(h)$ and $f_2(h)$. 
To this purpose we use that thanks to $4r\le 1$ we have
$\eta_{2r}\lesssim\eta\eta_{4r}$ for our exponential cut-off so that
$\int\eta_{2r} f_1 f_2 dh\lesssim\sup_{h}(\eta f_1)\int\eta_{4r} f_2 dh$.
We claim that for the first factor in the second term on the right hand side of (\ref{p114a}) we have
\begin{equation}\label{p116}
\sup_{h}\eta(h)\int_{-1}^0\int\eta(u^h-u)^2dxdt
\lesssim\int_{-1}^0\int\eta(u-\int_{-1}^0\int\eta u)^2dxdt
=D(u,1).
\end{equation}
Before inserting them, we give the easy arguments for (\ref{p115a}) and (\ref{p116}): By scaling we may assume
$r=1$ so that (\ref{p115a}) follows from Jensen's inequality in form of
\begin{equation}\nonumber
\int\eta(u-\int\eta u)^2dx\le \int\int\eta(x)\eta(x+h)(u^h(x)-u(x))^2dxdh
\end{equation}
and the fact that for our exponential cut-off $\eta(x)\eta(x+h)=\frac{1}{4}\exp(-(|x|+|x+h|))
\le\frac{1}{4}\exp(-\frac{1}{2}(|h|+|x|))=4\eta_2(h)\eta_2(x)$. For (\ref{p116}), by the triangle inequality
in $L^2$, it is enough to show for a constant $c$ ($\int_{-1}^0\int\eta udx dt$ in our case)
\begin{equation}
\sup_h\eta(h)\int\eta(u^h-c)^2dx\le \int\eta(u-c)^2dx.
\end{equation}
This inequality follows from writing
\begin{equation}\nonumber
\sup_h\eta(h)\int\eta(u^h-c)^2dx
=\sup_h\int\eta(h)\eta(x-h)(u(x)-c)^2dx
\end{equation}
and the fact that for our exponential cut-off $\eta(h)\eta(x-h)\le
\eta(x)$. Inserting (\ref{p115a}) and (\ref{p116}) into (\ref{p114a}) we obtain 
\begin{align}\label{p117}
\notag
&\av_{-r^2}^0\int\eta_{r}(u-\int\eta_r u)^2dxdt
 \lesssim  r^{ 2\alpha_0}  \Dprime  \\
& \qquad \Big(1+\frac{L^2}{r^3}D^2(u,1)\Big)
\int\eta_{4r}(h)\sup_{(t,x)\in(-1,0)\times\mathbb{R}}\eta(g^h-g)^2 dh.
\end{align}

\medskip

Before taking the (restricted) expectation of this inequality, we note that our choice 
of $c = \int_{-1}^0 \int \eta u dx dt$ depends on $u$ and so does our choice of 
coefficient $a_0 = \pi'(c)$.
Therefore $g$  has to be viewed as the solution of a stochastic heat equation with 
constant but \emph{random}, non-adapted 
coefficients and the standard regularity estimates do not apply immediately. This problem 
is addressed in Lemma~\ref{L3} where a bound on $\langle \sup_{a_0 \in [\lambda,1]} \sup_{(-1,0)\times\mathbb{R}} \eta(g^h-g)^2 \rangle_1$
is provided. So when taking the (restricted) expectation of (the square root of) (\ref{p117}) 
and inserting this estimate provided by Lemma \ref{L3}  we obtain
%
%
\begin{eqnarray}
\lefteqn{\Big\langle\Big(\av_{-r^2}^0\int\eta_{r}(u-\int\eta_r u)^2dxdt\Big)^\frac{1}{2}\Big\rangle_1}
\nonumber\\
&\lesssim&r^{\an} \Big( 1+ D'(u,1) +\frac{L}{r^\frac{3}{2}}\langle D^2(u,1)\rangle_1^\frac{1}{2}\Big)
 .
\end{eqnarray}
We now appeal to the triangle inequality in form of
\begin{eqnarray}
D(u,r)
&\le&\Big(\av_{-r^2}^0(\int\eta_{r}u-\av_{-r^2}^0\int\eta_r u)^2dt\Big)^\frac{1}{2}\nonumber\\
&&+\Big(\av_{-r^2}^0\int\eta_{r}(u-\int\eta_r u)^2dxdt\Big)^\frac{1}{2}\label{e:triangle}
\end{eqnarray}
and Lemma \ref{L4} for the upgrade to
\begin{eqnarray}
\langle D(u,r)\rangle_1&\lesssim&r^\frac{1}{2}+
r^{\an}\Big( 1 + D'(u,1) +\frac{L}{r^\frac{3}{2}}\langle D^2(u,1)\rangle_1^\frac{1}{2}\Big)
,\nonumber
\end{eqnarray}
which we rewrite as
\begin{eqnarray}\label{p119}
\langle D(u,r)\rangle_1&\lesssim&
r^{\an} \Big( 1 +D'(u,1) +\frac{L}{r^\frac{3}{2}}\langle D^2(u,1)\rangle_1^\frac{1}{2}\Big).
\end{eqnarray}
In this form, we see that (\ref{p119}) does not just hold for $r\ll 1$ but trivially
for $r\le1 $ with $r\sim 1$, since $D(u,r)\le\frac{1}{r^3}D(u,1)$. It remains to
appeal to  Lemma \ref{L2}.

\bigskip

{\sc Proof of Proposition \ref{P2}}.
For conciseness, we ignore the massive term in (\ref{i.1}) and fix $r\le 1$.
Following the argument in the proof of Lemma \ref{L0},
we first claim that the proposition reduces to the following deterministic estimate
\begin{equation}\label{p26}
D^2(\delta u,r)\lesssim \frac{1}{r^3}\int(\delta\xi)^2dxdt
\end{equation}
for any decaying $\delta u$ and $\delta\xi$ supported for $t\in(-1,0)$ related via
\begin{equation}\label{p211}
\partial_t\delta u-\partial_x^2(a\delta u)=\delta\xi,
\end{equation}
where $a:=\pi'(u)$ satisfies (\ref{p21}). Indeed, we note that by duality w.\ r.\ t.\ to 
the inner product $(g,f)\mapsto \av_{-r^2}^0\int\eta_r gfdxdt$,
\begin{eqnarray}\label{p210}
\lefteqn{D(u,r)=\sup\Big\{\;D(u,f):=\av_{-r^2}^0\int\eta_r u fdxdt\;\Big|}\\
\;&&\av_{-r^2}^0\int\eta_r f^2dxdt=1,\;\;{\rm supp}f\subset(-r^2,0)\times\mathbb{R},
\;\;\av_{-r^2}^0\int\eta_r fdxdt=0\;\Big\}.\nonumber
\end{eqnarray}
By the chain rule for the Malliavin derivative we thus obtain
\begin{equation}\nonumber
|\nabla D(u,r)|_1\le\sup_f|\nabla D(u,f)|_1,
\end{equation}
where the supremum runs over the set implicitly defined in (\ref{p210}), so that it
is enough to show for a fixed $f$
\begin{equation}\nonumber
|\nabla D(u,f)|_1^2\lesssim \frac{1}{r^3}.
\end{equation}
By definition (\ref{i.a}) of the carr\'e-du-champs of the Malliavin derivative 
applied to the linear functional $u\mapsto D(u,f)$, this amounts to show
\begin{equation}\nonumber
\av_{-r^2}^0\int\eta_r \delta u fdxdt
\lesssim \frac{1}{r^{\frac32}}\Big(\int(\delta\xi)^2dxdt\Big)^\frac{1}{2},
\end{equation}
where the infinitesimal perturbation $\delta u$ of the solution is related
to the infinitesimal perturbation $\delta\xi$ of the noise supported
on $(-1,0)\times\mathbb{R}$ via (\ref{p211}).
By the characterizing properties of the $f$'s, cf.\ (\ref{p210}), this estimate
in turn amounts to (\ref{p26}).

\medskip

In order to see \eqref{p26} we use the trivial estimate
\begin{align*}
D^2(\delta u,r)\ \lesssim \frac{1}{r^3} \int_{-1}^0 \int  \delta u(t,x)^2 dx dt 
\end{align*}
and apply Proposition~\ref{P4} for $h = \int_{-1}^{\cdot} \delta \xi ds$. Observing that 
$\int_{-1}^t$ is a bounded operator on $L^2(-1,0)$ we obtain \eqref{p26}.

\bigskip

{\sc Proof of Theorem \ref{T}}. 
In this proof $\lesssim$ and $\ll$ refer to constants only depending on $\lambda$ and eventually on $L$ and $\alpha$.
For some  $\theta \in (0,1)$ to be chosen later  we consider the random variable $D(u, \theta)$. By Proposition~\ref{P2}
we know that $D(u,\theta)$ is a Lipschitz variable with respect to perturbations of  the noise 
which are supported in $(-1,0) \times \mathbb{R}$ and we have
\begin{equation}\nonumber
|\nabla D(u,\theta)|_1\lesssim \theta^{-\frac32}.
\end{equation}

By concentration of measure, cf.\ the beginning of Section \ref{S2}, applied to the restricted ensemble $\langle \cdot \rangle_1$ we
 conclude that suitably rescaled fluctuations
 \begin{equation*}
 \chi =\theta^{\frac32}\big( D(u, \theta) - \langle D(u,\theta)\rangle_1 \big)
 \end{equation*}
satisfy
 $\langle\exp(\frac{1}{C}\chi^2)\rangle_1\le 2$,
and thus a fortiori $\langle\exp(\frac{1}{C}\chi^2)\rangle\le 2$.
Combining this with  Proposition \ref{P1} we get the almost-sure inequality
\begin{align*}
\frac{1}{\theta^{\an}} D(u,\theta) \lesssim 1+ D'(u,1) +\frac{L}{\theta^\frac{3}{2}}(D'(u,1)+1)+ \frac{1}{\theta^{\frac32 +\an}} \chi .
\end{align*}

By the invariance in law under the scaling (\ref{t1}) \& (\ref{t2}) \& (\ref{i.b}),
this yields for any length scale $R \leq 1$ 
\begin{align}\label{t3a}
\frac{1}{\theta^{\an}}  \frac{1}{R^{\frac12}}D(u, \theta R ) 
	\lesssim 1+ \frac{1}{R^{\frac12}} D'(u,R) +\frac{R^{\frac12}L}{\theta^\frac{3}{2}} \big(\frac{1}{R^{\frac12}}D'(u,R)+1\big)
		+ \frac{1}{\theta^{\frac32 +\an}} \chi_R ,
\end{align}
with an $R$-dependent random variable $\chi_R$ of Gaussian moments 
$\langle\exp(\frac{1}{C}\chi_R^2)\rangle\lesssim 1$. 
Then using the fact that 
$\Dprime(u,R)=\int\eta_R(u-\int\eta_Ru)^2dx_{|t=-R^2}$ satisfies 
\begin{align}
\av_\frac{R}{2}^RD'(u,R')dR' & = \av_{\frac{R}{2}}^R \Big( \int\eta_{R'}(u-\int\eta_{R'}u)^2dx_{|t=-R^2} \Big)^{\frac12} dR' \notag\\
&\leq\frac{2}{R} \av_{\frac{R}{2}}^R \Big( \int\eta_{R'}(u-\av_{-R}^0\int \eta_{R}u)^2dx_{|t=-R^2} \Big)^{\frac12} R' dR' \notag\\
&\lesssim \av_{-R^2}^0  \Big( \int\eta_{R}(u-\av_{-R}^0\int \eta_{R}u)^2dx  \Big)^{\frac12}dt \lesssim D(u,R).\label{dprimevsd}
\end{align}
 we see that by replacing $R$ by $R'$ in (\ref{t3a}) and by
averaging over $R'\in(\frac{R}{2},R)$ we obtain
\begin{align*}
\frac{1}{\theta^{\an}} \frac{1}{R^{\frac12}}  D(u, \theta R ) 
	\lesssim 1+ \frac{1}{R^{\frac12}} D(u,R)+\frac{L}{\theta^\frac{3}{2}}( D(u,R)+R^{\frac12})  
		+ \frac{1}{\theta^{\frac32 +\an}} \chi_R' ,
\end{align*}
where $\chi'_R:=\av_\frac{R}{2}^R\chi_{R'}dR'$ still has Gaussian moments
$\langle\exp(\frac{1}{C}{\chi'_R}^2)\rangle\le 2$, since the latter property is preserved
by convex combination. In order to prepare for recursion we introduce an $\alpha<\an$ and 
rewrite the last estimate as
\begin{align*}
\frac{1}{(\theta R)^{\alpha }} D(u,  \theta R) 
	\lesssim & (R^{\frac12 - \alpha} \theta^{\an - \alpha}   + \frac{LR^{1-\alpha}}{\theta^{\frac32 -(\an - \alpha)}}  ) 
	+ (\theta^{\an-\alpha}  + \frac{R^{\frac{1}{2}} L}{\theta^{\frac{3}{2}- (\an- \alpha) } } ) \frac{1}{R^{\alpha } }  D(u,R) \\
		&+ \frac{R^{\frac{1}{2} -\alpha}}{\theta^{\frac{3}{2}+ \alpha }} \chi_R'  .
\end{align*}
Thanks to $\alpha < \alpha_0$ we can first fix $\theta \ll 1$ and then $R_0 \ll1$ so that this estimate turns into
\begin{align}\label{T455}
\frac{1}{(\theta R)^{\alpha}} D(u,  \theta R) 
	\leq  \frac12 \frac{1}{R^{\alpha} }  D(u,R) 
		+  C R^{\frac12  - \alpha } \chi_R'+C 
\end{align}
for all $R\leq  R_0$.
This prompts to consider the random variable 
\begin{equation}\nonumber
\bar\chi_{R_0}:=\max_{n=0,1,\cdots}(\theta^n)^{\frac{1}{2}-\alpha}\chi'_{R_0 \theta^n},
\end{equation}
%
which in view of (recall $\alpha<\frac{1}{2}$)
\begin{eqnarray}\nonumber
\bar\chi_{R_0}&\le&
\sum_{n=0}^\infty(\theta^n)^{\frac{1}{2}-\alpha}\chi'_{\theta^n R_0}\\
&=&\frac{1}{1-\theta^{\frac{1}{2}-\alpha}}\times\mbox{convex combination of}
\;\{\chi'_{\theta^n R_0}\}_{n=0,1,\cdots}\label{conv-comb}
\end{eqnarray}
has Gaussian moments $\big\langle\exp\big(\frac{1}{C}\bar\chi^2\big)\big\rangle\leq 2$
since by construction, the random variables $\{\chi'_R\}_{R}$ have a uniform Gaussian moment bounds.
Thanks to the factor $\frac{1}{2}<1$ the estimate \eqref{T455} can be iterated to yield
\begin{eqnarray*}
\sup_{n=0,1,\cdots}\frac{1}{(\theta^n )^{\alpha}}D(u,R_0 \theta^n )
\lesssim  D(u,R_0)+ \bar\chi +1,
\end{eqnarray*}
which implies (using once more $D(u,r)\le (\frac{R}{r})^\frac{3}{2}D(u,R)$ for any scales $r\le R$
to bridge the dyadic gaps as well as the gap between $R_0$ and $1$)
\begin{eqnarray*}
\sup_{r\le1}\frac{1}{r^{\alpha}}D(u,r)
\lesssim  D(u,1)+ \bar\chi +1.
\end{eqnarray*}
Theorem~\ref{T} now follows by invoking Lemma~\ref{L0} and relabelling $\alpha$ as $\alpha_0$.

\bigskip

{\sc Proof of Proposition \ref{P3}}.
We start by observing
\begin{eqnarray}\nonumber
\frac{d}{dt}\frac{1}{2}\int\eta (\partial_xv)^2dx
&=&-\int\partial_x(\eta\partial_xv)(a_0\partial_x^2v+f)dx\nonumber\\
&=&-\int \eta(a_0(\partial_x^2v)^2+f\partial_x^2v)dx
-\int\partial_x\eta\partial_xv(a_0\partial_x^2v+f)dx,\nonumber
\end{eqnarray}
so that because of $a_0\in[\lambda,1]$ and $|\partial_x\eta|\le\eta$
we obtain by Young's inequality
\begin{eqnarray}\label{p32}
\frac{d}{dt}\frac{1}{2}\int\eta (\partial_xv)^2dx
&\le&-\int\eta(\lambda(\partial_x^2v)^2+f\partial_x^2v)dx
+\int\eta|\partial_xv|(|\partial_x^2v|+|f|)dx\nonumber\\
&\le&-\frac{1}{C}\int\eta(\partial_x^2v)^2dx+C\int\eta((\partial_xv)^2+f^2)dx.
\end{eqnarray}
%
%
%
Dropping the good r.\ h.\ s.\ term, we rewrite this as
\begin{eqnarray}\nonumber
\frac{d}{dt}(t+1)\int\eta (\partial_xv)^2dx
&\lesssim&\int\eta((\partial_xv)^2+f^2)dx,
\end{eqnarray}
so that we obtain from integration in $t\in(-1,0)$
\begin{equation}\label{p46}
\sup_{t\in(-1,0)}(t+1)\int\eta (\partial_xv)^2dx\lesssim\int_{-1}^0\int\eta((\partial_xv)^2+f^2)dxdt.
\end{equation}

\medskip

Thanks to the constant coefficients, also the (localized) $L^2$-norm is well-behaved.
Indeed, from (\ref{i.6}) we obtain
\begin{eqnarray}\nonumber
\frac{d}{dt}\frac{1}{2}\int\eta v^2dx
&=&\int\eta v(a_0\partial_x^2v+f)dx\nonumber\\
&=&\int \eta(-a_0(\partial_xv)^2+vf)dx
-a_0\int\partial_x\eta v\partial_xvdx,\nonumber
\end{eqnarray}
so that because of $a_0\in[\lambda,1]$ and $|\partial_x\eta|\le\eta$, we obtain by Young's inequality
\begin{eqnarray}\nonumber
\frac{d}{dt}\int\eta v^2dx
&\le&-\frac{1}{C}\int \eta(\partial_xv)^2dx
+\int\eta(v^2+f^2)dx.\nonumber
\end{eqnarray}
From the integration in $t$ of this differential inequality for $\int\eta v^2dx$ we learn
\begin{equation}\label{p39}
\sup_{t\in(-1,0)}\int\eta v^2dx+\int_{-1}^0\int\eta(\partial_xv)^2dxdt
\lesssim\int_{-1}^0\int\eta f^2dxdt+\int\eta v^2dx_{|t=-1}.
\end{equation}
The combination of this with (\ref{p46}) yields
\begin{equation}\label{p36}
\sup_{t\in(-1,0)}\big((1+t)\int\eta(\partial_xv)^2+\int\eta v^2dx\big)
\lesssim\int_{-1}^0\int\eta f^2dxdt+\int\eta v^2dx_{|t=-1}.
\end{equation}

\medskip

In view of this a priori estimate (\ref{p36}),
the desired estimate (\ref{p35})  follows from the embedding
\begin{equation}\nonumber
\sup_{(t,x)\in(-1,0)\times\mathbb{R}}(1+t)^\frac{1}{2}\eta v^2
\lesssim \sup_{t\in(-1,0)}\big((1+t)\int\eta(\partial_xv)^2dx+\int\eta v^2dx\big),
\end{equation}
which easily is is seen to hold: Because of
\begin{equation}\nonumber
\sup_x\eta v^2\lesssim\int|\partial_x(\eta v^2)|dx\lesssim
\int(\eta|v\partial_xv|+|\partial_x\eta|v^2)dx
\lesssim\int\eta(|\partial_xv|+ |v|)|v|dx,
\end{equation}
we obtain by Young's inequality for $t\in(-1,0)$
\begin{eqnarray*}\nonumber
(t+1)^\frac{1}{2}\sup_x\eta v^2
&\lesssim&(1+t)\int\eta((\partial_xv)^2+ v^2)dx+\int\eta v^2dx\nonumber\\
&\lesssim&(1+t)\int\eta(\partial_xv)^2dx+\int\eta v^2dx.
\end{eqnarray*}

\bigskip

\medskip 
{\sc Proof of Proposition \ref{P1old}}.
As before we ignore the massive term in (\ref{i.1}) and monitor $\delta u:=u^h-u$, where $u^h(t,x)=u(t,x+h)$ denotes
a spatial shift of the stationary solution of (\ref{i.1}). As in the proof of Proposition~\ref{P1}  $\delta u$
satisfies the formally linear equation \eqref{p110}.

\medskip

We will again derive a recursive estimate to go from scale one to scale $r$ for the 
functions $\delta u$. This time it is given by the estimate
\begin{eqnarray}
\av_{-r^2}^0\int\eta_r(\delta u)^2dxdt
&\lesssim&\Big(1+\frac{1}{r^3}\int_{-1}^0\int\eta(a_h-a_0)^2dxdt\Big)\nonumber\\
&\times&\Big(\int_{-1}^0\int\eta(\delta u)^2dxdt
+\sup_{(t,x)\in(-1,0)\times\mathbb{R}}\eta(\delta g)^2\Big),\label{p12old}
\end{eqnarray}
which we shall establish for all $r\ll 1$. Compared to \eqref{p12} the term 
$$r^{ 2\an} \big(\frac{h}{r}\big)^2e^{|h|} \Dprime$$ on the right hand side, which limits that 
H\"older regularity to the strictly positive but possibly small $\an>0$ coming from Proposition ~\ref{P5},
has disappeared. The price to pay is the extra term 
$\int_{-1}^0\int\eta(\delta u)^2dxdt$ which will eventually lead to an  estimate 
which is quadratic in $u$.

\medskip

%
%
We observe that it is enough to establish for any $R\in[\frac{1}{2},1]$ the estimate
\begin{eqnarray*}
\av_{-r^2}^0\int\eta_r(\delta u)^2dxdt
&\lesssim&\Big(1+\frac{1}{r^3}\int_{-R^2}^0\int(t+R^2)^{-\frac{1}{2}}\eta(a_h-a_0)^2dxdt\Big)\\
&\times&\Big(
\int\eta(\delta u)^2dx_{|t=-R^2}
+\sup_{(t,x)\in(-R^2,0)\times\mathbb{R}}\eta(\delta g)^2\Big),\nonumber
\end{eqnarray*}
since the  integral of this estimate over $R\in[\frac{1}{2},1]$ yields
(\ref{p12}), using the integrability of $(t+R^2)^{-\frac{1}{2}}$ thanks to $\frac{1}{2}<1$.
To simplify notation, we replace $R\sim 1$ by unity, so that it remains to show
\begin{eqnarray}
\av_{-r^2}^0\int\eta_r(\delta u)^2dxdt
&\lesssim&\Big(1+\frac{1}{r^3}\int_{-1}^0\int(t+1)^{-\frac{1}{2}}\eta(a_h-a_0)^2dxdt\Big)\nonumber\\
&\times&\Big(
\int\eta(\delta u)^2dx_{|t=-1}
+\sup_{(t,x)\in(-1,0)\times\mathbb{R}}\eta(\delta g)^2\Big).\label{p16old}
\end{eqnarray}

\medskip

As before we split the 
solution $\delta u=\delta g+v+w$, but this time $v$ is defined through the 
constant-coefficient initial value problem
\begin{equation}\nonumber
\partial_tv-a_0\partial_x^2v=0\;\;\mbox{for}\;t\in(-1,0),\quad
v=\delta u\;\;\mbox{for}\;t=-1,
\end{equation}
%
%
and $w$ is defined through the initial value problem
\begin{equation}\label{p18old}
\partial_tw-\partial_x^2(a_hw)=\partial_x^2((a_h-a_0)(\delta g+v))\;\;\mbox{for}\;t\in(-1,0),\quad
w=0\;\;\mbox{for}\;t=-1.
\end{equation}
%

\medskip

In view  of the constant coefficient $a_0$, this time we can bound $v$ using  Proposition \ref{P3} rather than Proposition~\ref{P5}.
More precisely, the first bound \eqref{p35} of Proposition \ref{P3} implies that
\begin{equation}\label{p13}
\sup_{(t,x)\in(-1,0)\times\mathbb{R}}(1+t)^\frac{1}{2}\eta v^2\lesssim
\int\eta(\delta u)^2dx_{|t=-1},
\end{equation}
which implies in particular for $r\le\frac{1}{2}$ (which amounts to $r \leq \frac14$ before setting setting $R=1$
above)
\begin{equation}\label{p15}
\av_{-r^2}^0\int\eta_r v^2\lesssim
\int\eta(\delta u)^2dx_{|t=-1}.
\end{equation}
To bound $w$ we  use the first part (\ref{p43}) of Proposition \ref{P4} (with $\eta$ replaced by $\eta_\frac{1}{2}$) 
to gather that
\begin{equation}\nonumber
\int_{-1}^0\int\eta_{\frac{1}{2}}w^2dxdt\lesssim\int_{-1}^0\int\eta_\frac{1}{2}(a_h-a_0)^2((\delta g)^2+v^2)dxdt,
\end{equation}
which implies for $r\le\frac{1}{2}$ (using the obvious inequality $\av_{-r^2}^0\int\eta_r\cdot dxdt
\le(\frac{R}{r})^3\av_{-R^2}^0\int\eta_R\cdot dxdt$ for $r\le R$ and since $\eta_\frac{1}{2}\lesssim\eta^2$)
\begin{eqnarray}\label{p14}
\lefteqn{r^3\av_{-r^2}^0\int\eta_rw^2dxdt}\\
&\lesssim&
\int_{-1}^0(1+t)^{-\frac{1}{2}}\int\eta(a_h-a_0)^2dxdt
\sup_{(t,x)\in(-1,0)\times\mathbb{R}}(1+t)^\frac{1}{2}\eta ((\delta g)^2+v^2).\nonumber
\end{eqnarray}
Inserting (\ref{p13}) into (\ref{p14}) yields 
\begin{eqnarray}
\av_{-r^2}^0\int\eta_rw^2dxdt&\lesssim&
\frac{1}{r^3}\int_{-1}^0(1+t)^{-\frac{1}{2}}\int\eta(a_h-a_0)^2dxdt\nonumber\\
&\times&\Big(
\int\eta(\delta u)^2dx_{|t=-1}
+\sup_{(t,x)\in(-1,0)\times\mathbb{R}}\eta (\delta g)^2\Big).\nonumber
\end{eqnarray}
Finally, as in Proposition~\ref{P1} we have 
\begin{equation}\nonumber
\av_{-r^2}^0\int\eta_r(\delta g)^2dxdt\lesssim\sup_{(t,x)\in(-1,0)\times\mathbb{R}}\eta(\delta g)^2.
\end{equation}
Combining the two last estimates with (\ref{p15}) yields (\ref{p16old}) for $\delta u=\delta g+v+w$.

\medskip

The post-processing of (\ref{p12old}) follows the same lines as the corresponding argument in the proof of Proposition~\ref{P1}
and we only give a sketch. Making the choice of
$a_0=\pi'(c)$ with $c:=\int_{-1}^0\int\eta udxdt$ and using \eqref{e:ah-vs-anull}, the estimate \eqref{p12old} turns into
\begin{eqnarray}
\lefteqn{\av_{-r^2}^0\int\eta_{2r}(u^h-u)^2dxdt}\nonumber\\
&\lesssim&\Big(1+\frac{L^2}{r^3}\int_{-1}^0\int\eta(\delta u)^2dxdt+\frac{L^2}{r^3}D^2(u,1)\Big)\nonumber\\
&\times&\Big(\int_{-\frac{1}{2}}^0\int\eta(\delta u)^2dxdt
+\sup_{(t,x)\in(-1,0)\times\mathbb{R}}\eta(\delta g)^2\Big).\nonumber
\end{eqnarray}
Then we  integrate this estimat in $h$ according to $\int\eta_{2r}(h)\cdot dh$, and use \eqref{p115a} to compare the left hand side to $\av_{-r^2}^0\int\eta_{r}(u-\int\eta_r u)^2dxdt$
and \eqref{p116} to compare the supremum of first factor on the right hand side weighted by $\eta(h)$ to $1 + \frac{L^2}{r^3} D^2(u,1)$. In this 
way we obtain
\begin{eqnarray}\label{p117old}
\lefteqn{\av_{-r^2}^0\int\eta_{r}(u-\int\eta_r u)^2dxdt
\lesssim\Big(1+\frac{L^2}{r^3}D^2(u,1)\Big)}\nonumber\\
&\times&\int\eta_{4r}(h)\big(\int_{-\frac{1}{2}}^0\int\eta(\delta u)^2dxdt
+\sup_{(t,x)\in(-1,0)\times\mathbb{R}}\eta(\delta g)^2\big)dh.
\end{eqnarray}

\medskip

Then,  we take the restricted expectation of the square root of (\ref{p117old}). As 
in the proof of Proposition~\ref{P1} (keeping in mind the discussion 
about the dependence of $g$ on the random coefficient $a_0$), Lemma~\ref{L3} allows us to control the term involving $g$ to obtain a factor $r^\alpha$ for any 
$\alpha < \frac12$. For the term involving $\delta u$ in the second factor we invoke Lemma~\ref{L1} and obtain
\begin{eqnarray}
\lefteqn{\Big\langle\Big(\av_{-r^2}^0\int\eta_{r}(u-\int\eta_r u)^2dxdt\Big)^\frac{1}{2}\Big\rangle_1}
\nonumber\\
&\lesssim&\Big(1+\frac{L}{r^\frac{3}{2}}\langle D^2(u,1)\rangle_1^\frac{1}{2}\Big)
\Big(r^\frac{1}{2}+r\langle D^2(u,1)\rangle_1^\frac{1}{2}+r^\alpha\Big).\nonumber
\end{eqnarray}
Then as in the proof of Proposition~\ref{P1} we make use of the triangle inequality~\eqref{e:triangle} 
and Lemma \ref{L4} to replace the term $\Big\langle\Big(\av_{-r^2}^0\int\eta_{r}(u-\int\eta_r u)^2dxdt\Big)^\frac{1}{2}\Big\rangle_1$ on the left hand side by $\langle D(u,r)\rangle_1$. In this way we 
finally obtain 
\begin{eqnarray*}
\langle D(u,r)\rangle_1&\lesssim&
r^\alpha\Big(1+\frac{L}{r^\frac{3}{2}}\langle D^2(u,1)\rangle_1^\frac{1}{2}\Big)
\Big(1+r^\frac{1}{2}\langle D^2(u,1)\rangle_1^\frac{1}{2}\Big)
\end{eqnarray*}
which yields the desired estimate after
appealing to  Lemma \ref{L2}.

\bigskip

\bigskip

{\sc Proof of Theorem \ref{T2}}. 
In this proof $\lesssim$ and $\ll$ refer to generic constants only depending on $\lambda$,
$\alpha$,  and eventually $L$. The first steps in this proof resemble the argument used in the proof of Theorem~\ref{T}:
According to Proposition \ref{P1} we have
\begin{equation}\nonumber
\frac{1}{\theta^\alpha}\langle  D(u,\theta)\rangle_1\lesssim\big(1+\frac{L}{\theta^\frac{3}{2}}(1+D'(u,1))\big)\,
\big(1+\theta^\frac{1}{2}D'(u,1)\big)
\end{equation}
and according to Proposition \ref{P2}
\begin{equation}\nonumber
|\nabla  D(u,\theta)|_1\lesssim \theta^{-\frac32}.
\end{equation}
As in the proof of Theorem~\ref{T}, we apply the concentration of measure
 to the restricted ensemble $\langle\cdot\rangle_1$
and obtain the existence
of a random variable $\chi$ with Gaussian bounds,  $\langle\exp(\frac{1}{C}\chi^2)\rangle\le 2$ such that
\begin{equation}\nonumber
\frac{1}{\theta^\alpha}D(u,\theta)
\lesssim \big(1+\frac{L}{\theta^\frac{3}{2}}(1+D'(u,1))\big)\big(1+\theta^\frac{1}{2}D'(u,1)\big) +\frac{1}{\theta^{\frac32+\alpha}} \chi,
\end{equation}
where  we think of small $\theta \ll 1$. 
We apply the invariance in law under the scaling (\ref{t1}) \& (\ref{t2}) \& (\ref{i.b}),
to obtain for any length scale $R$
\begin{align}\label{t3}
 \frac{1}{\theta^\alpha}\frac{1}{R^\frac{1}{2}}D(u,\theta R)  
\lesssim \big(1+\frac{1}{\theta^\frac{3}{2}}(R^\frac{1}{2}+D'(u,R))\big)
\big(1+\frac{\theta^\frac{1}{2}}{R^\frac{1}{2}}D'(u,R)\big)+ \frac{1}{\theta^{\frac32+\alpha}}\chi_R,
\end{align}
with an $R$-dependent random variable $\chi_R$ of Gaussian moments 
$\langle\exp(\frac{1}{C}\chi_R^2)\rangle\leq 2$. As in the proof of Theorem~\ref{T} we integrate out the 
initial time $-R$ in this estimate using \eqref{dprimevsd}  to get $\av_\frac{R}{2}^RD'(u,R')dR'
\lesssim D(u,R)$ 
and get
\begin{align*}
  \frac{1}{\theta^\alpha}\frac{1}{R^\frac{1}{2}}D(u,\theta R) 
\lesssim \big(1+\frac{1}{\theta^\frac{3}{2}}(R^\frac{1}{2}+D(u,R))\big)
\big(1+\frac{\theta^\frac{1}{2}}{R^\frac{1}{2}}D(u,R)\big)+\frac{1}{\theta^{\frac32+\alpha}}\chi'_R,
\end{align*}
where $\chi'_R:=\av_\frac{R}{2}^R\chi_{R'}dR $ still has Gaussian moments
$\langle\exp(\frac{1}{C}{\chi'_R}^2)\rangle\le 2$.  
\medskip

At this point our argument deviates from the proof of Theorem~\ref{T}, because the bound on the 
right hand side is non-linear in $D(u, R)$.
We reformulate it as
\begin{equation}\nonumber
\left.\begin{array}{c}D(u,R)\le \theta^{\frac32}\\
R\le \theta^{3} \end{array}\right\}
\;\;\Longrightarrow\;\;
\frac{1}{\theta^\alpha}\frac{1}{R^\frac{1}{2}}D(u,\theta R)
\le C\frac{\theta^\frac{1}{2}}{R^\frac{1}{2}}D(u,R) +\frac{1}{\theta^{\frac32+\alpha}}\chi''_R,
\end{equation}
where $\chi''_R\sim\chi'_{R}+1$ still has Gaussian moments
$\langle\exp(\frac{1}{C}{\chi''_R}^2)\rangle\leq 2$. Hence selecting $\theta\sim 1$
sufficiently small, we obtain
\begin{equation}\label{t5}
\left.\begin{array}{c}D(u,R)\le \theta^{\frac32}\\R\le \theta^3\end{array}\right\}
\;\;\Longrightarrow\;\;
\frac{1}{(\theta R)^\alpha}D(u,\theta R)
\le\frac{1}{2}\frac{1}{R^\alpha}D(u,R) + \frac{R^{\frac{1}{2}-\alpha}}{\theta^{\frac32+\alpha}} \chi''_R.
\end{equation}
Since (\ref{t5}) implies in particular 
$D(u,\theta R)\le \frac{R^\frac{1}{2}}{\theta^{\frac32}}\chi''_R+\frac{1}{2} D(u,R)$,
we see that in order to convert (\ref{t5}) into a self-propelling iteration, 
we need $R^\frac{1}{2}\chi''_R\le\frac{1}{2}\theta^3$,
$(\theta R)^\frac{1}{2}\chi''_{\theta R}\le\frac{1}{2}\theta^3$ and so on.
Now we can proceed as in the proof of Theorem~\ref{T} and consider the random variable 
\begin{equation}\nonumber
\bar\chi_R:=\max_{n=0,1,\cdots}(\theta^n)^{\frac{1}{2}-\alpha}\chi''_{\theta^n R}
\ge\max_{n=0,1,\cdots}(\theta^n)^{\frac{1}{2}}\chi''_{\theta^n R},
\end{equation}
which thanks to $\alpha'< \frac12$ can be written as a multiple of a convex combination of  $\chi_{\theta^n R}$ as in \eqref{conv-comb},
and which therefore 
has Gaussian moments $\big\langle\exp\big(\frac{1}{C}\bar\chi^2\big)\big\rangle\leq 2$ itself.
From (\ref{t5}) we learn 
\begin{eqnarray*}
\lefteqn{D(u,R)\le \theta^\frac{3}{2},\quad R \leq \theta^3 \quad\mbox{and}\quad R^\frac{1}{2}\bar\chi_R\le\frac{1}{2} \theta^3}\\
&\Longrightarrow&\forall\;n\in\mathbb{N}\;\;
\frac{1}{(\theta^n R)^\alpha}D(u,\theta^n R)
\le \frac{1}{2}\frac{1}{(\theta^{n-1}R)^\alpha}D(u,\theta^{n-1}R)
+\frac{R^{\frac{1}{2}-\alpha}}{\theta^{\frac32+\alpha}}\bar\chi_{R}.
\end{eqnarray*}
Thanks to the factor $\frac{1}{2}<1$ the last statement can be iterated to yield
\begin{eqnarray*}
\lefteqn{D(u,R)\le \theta^\frac32, \quad R \leq \theta^3 \quad\mbox{and}\quad R^\frac{1}{2}\bar\chi_R\le \frac{1}{2} \theta^3}\\
&\Longrightarrow&
\sup_{n=0,1,\cdots}\frac{1}{(\theta^n R)^\alpha}D(u,\theta^n R)
\le  \frac{1}{2}\frac{1}{R^\alpha}D(u,R)+ 2\frac{R^{\frac{1}{2}-\alpha}}{\theta^{\frac32+\alpha}}\bar\chi_{R},
\end{eqnarray*}
which implies (using once more $D(u,r)\le (\frac{R}{r})^\frac{3}{2}D(u,R)$ for any scales $r\le R$
to bridge the dyadic gaps)
\begin{align*}
&D(u,R)\le \theta^{\frac{3}{2}}, \quad R \leq \theta^3 \quad\mbox{and}\quad R^\frac{1}{2}\bar\chi_R\le\frac{1}{2}\theta^3\\
& \quad \Longrightarrow \quad 
\sup_{r\le R}\frac{1}{r^\alpha}D(u,r)
\le \frac{1}{2R^\alpha} + \frac{1}{(R\theta)^\alpha}.
\end{align*}

\medskip

Summing up, we learned that for any length scale $R \ll 1$,  we have for some constant $C_0=C_0(\lambda,\alpha,L)$ whose value we want to 
momentarily remember
\begin{eqnarray}\label{e:T21}
D(u,R)\leq \frac{1}{C_0}\quad\mbox{and}\quad R^\frac{1}{2}\bar\chi_{R} \leq \frac{1}{C_0}
&\Longrightarrow&
\sup_{r\le R}\Big(\frac{R}{r}\Big)^\alpha D(u,\rho)
\lesssim 1.
\end{eqnarray}
We now apply the elementary inequality  
\begin{align*}
\sup_{r \leq 1} \frac{1}{r^\alpha}D(u,r) \leq  \frac{1}{R^\alpha}\sup_{0 \leq r \leq R} \Big(\frac{R}{r}\Big)^\alpha D(u,
r) + \frac{1}{R^{\alpha - \alpha_0}} \sup_{0 \leq r \leq 1} \frac{1}{r^{\alpha_0}} D(u,r), 
\end{align*}
which is valid for any $\alpha_0 \leq \alpha$ and which we will use for an exponent $\alpha_0$ appearing in Theorem~\ref{T}.
%
As a value for $R$ we choose $R =\frac{1}{\bar{M}}$
where
\begin{align*}
\bar{M} &:= \inf \Big\{\frac{1}{r}   \Big|  D(u, r) \leq \frac{1}{C_0}  \quad \text{and} \quad  r^{\frac12 }\chi_r \leq \frac{1}{C_0  } \Big\}\\
&\leq \max \Big\{  \Big( C_0 \sup_{r \leq 1} \frac{1}{r^{\alpha_0}} D(u,r)\Big)^{\frac{1}{\alpha_0}} , \inf \big\{  M>0 \big| \bar{\chi}_{\frac{1}{M}} \leq \frac{M^{\frac12}}{C_0}   \big\} \Big\}.
 \end{align*}
With these choices \eqref{e:T21} turns into
\begin{align}
\sup_{r \leq 1} \frac{1}{r^\alpha}D(u,r) &\lesssim  \bar{M}^\alpha  + \bar{M}^{\alpha - \alpha_0} \sup_{0 \leq r \leq 1} \frac{1}{r^{\alpha_0}} D(u,r)\notag\\
&\lesssim  \bar{M}^\alpha  + \Big(\sup_{0 \leq r \leq 1} \frac{1}{r^{\alpha_0}} D(u,r) \Big)^{\frac{\alpha}{\alpha_0}},\label{T22}
\end{align} 
where in the last step we have used Young's inequality.
We claim that the random variable $\bar{M}$ has stretched exponential moments of the form 
\begin{align}
\big\langle \exp\big(\frac{1}{C} \bar{M}^{2 \alpha_0} \big)  \big\rangle \leq 2. \label{M11}
\end{align} 
Indeed, on the one hand
\begin{align}\label{M22}
\exp \Big( \frac{1}{C}\Big( C_0 \sup_{r \leq 1} \frac{1}{r^{\alpha_0}} D(u,r)\Big)^{\frac{2 \alpha_0}{\alpha_0}} \Big) \leq 2
\end{align} 
 by Theorem~\ref{T}, and on the other hand we have for any threshold $\mu$
 \begin{align*}
\Big\langle I\Big( \inf \Big\{  M  >0 \colon  \bar{\chi}_{\frac{1}{M}} \leq \frac{M^{\frac12}}{C_0} \Big\} \geq \mu \Big) \Big\rangle &\leq \Big\langle I\Big( \bar{\chi}_{\frac{1}{\mu}} \geq \frac{ \mu^{\frac{1}{2} }}{C_0} \Big) \Big\rangle \lesssim e^{-\frac{\mu}{C}},
\end{align*}
due to the uniform Gaussian tails of the $\bar{\chi}_{\frac{1}{\mu}}$, which implies that
\begin{align*}
\exp \Big(\frac{1}{C} \inf \Big\{  M  >0 \colon  \bar{\chi}_{\frac{1}{M}} \leq \frac{M^{\frac12}}{C_0} \Big\}\Big) \leq 2.
\end{align*}
Using \eqref{M11} and \eqref{M22} on \eqref{T22} yields as desired
\begin{align*}
\Big\langle \exp \Big( \frac{1}{C}\Big(  \sup_{0 \leq f \leq 1} \frac{1}{r^{\alpha}} D(u,r) \Big)^{2\frac{\alpha_0}{\alpha}} \Big)\Big\rangle \leq 2.
\end{align*}

%

\bigskip

\bigskip

{\sc Proof of Lemma \ref{L1}}.
Since thanks to $r\ll 1$ we have $\int\eta_r(h)(e^{|h|}-1)^2dh\lesssim r^2$, it is enough to show
for any shift $h$
\begin{eqnarray}\label{l11}
\lefteqn{
\Big\langle\int_{-\frac{1}{2}}^0\int\tilde\eta^2(u^h-u)^2dxdt\Big\rangle_1
}\nonumber\\
&\lesssim&|h|+(e^{|h|}-1)^2\Big\langle\int_{-1}^0\int\tilde\eta^2(u-c)^2dxdt\Big\rangle_1,
\end{eqnarray}
where, as in Lemma \ref{L0}, 
for the upcoming calculations we have replaced the exponential cut-off $\eta=\eta_2^2$ by its smooth version
$\tilde\eta^2$ where
\begin{equation}\label{l21}
\tilde\eta(x):=\exp(-\frac{1}{2}\sqrt{x^2+1})\sim\eta_2(x)
\end{equation}
and we have set for abbreviation $c:=\int_{-1}^0\int\tilde\eta^2 udxdt$. By the martingale argument
based on the stochastic differential equation
\begin{equation}\nonumber
\partial_t(u^h-u)=-(-\partial_x^2)(\pi(u^h)-\pi(u))+(\xi^h-\xi)
\end{equation}
we have
\begin{eqnarray}\label{l12}
\lefteqn{\frac{d}{dt}\frac{1}{2}\Big\langle\int\tilde\eta(u^h-u)(1-\partial_x^2)^{-1}\eta(u^h-u)dx\Big\rangle_1}
\nonumber\\
&=&-\Big\langle\int\tilde\eta(u^h-u)(1-\partial_x^2)^{-1}\tilde\eta(-\partial_x^2)(\pi(u^h)-\pi(u))dx
\Big\rangle_1\nonumber\\
&&+\frac{1}{2}\int\tilde\eta(\tilde\eta-\frac{1}{2}\tilde\eta^he^{-|h|}-\frac{1}{2}
\tilde\eta^{-h}e^{-|h|})dx.
\end{eqnarray}
Let us make two comments on (\ref{l12}): Like in Lemma \ref{L0} we use physics notation
in the sense that an operator acts on all the terms to its right, e.\ g.\
in the above expression
$(1-\partial_x^2)^{-1}\tilde\eta(-\partial_x^2)(\pi(u^h)-\pi(u))=
(1-\partial_x^2)^{-1}[\tilde\eta(-\partial_x^2)(\pi(u^h)-\pi(u))]$.
The last term in (\ref{l12}), which comes from the quadratic variation
of the white noise in time, cf.\ Lemma \ref{L0} for a heuristic discussion, assumes this form because 
$\frac{1}{2}\exp(-|x-y|)$ is the (translation-invariant) kernel of the operator $(1-\partial_x^2)^{-1}$, 
so that $\tilde\eta(x)\frac{1}{2}\exp(-|x-y|)\tilde\eta(y)$ is the kernel of the operator 
$\tilde\eta(1-\partial_x^2)^{-1}\tilde\eta$, so that the quadratic variation is indeed given by
\begin{eqnarray*}\nonumber
\lefteqn{\frac{1}{2}\int\int\tilde\eta(x)\frac{1}{2}\exp(-|x-y|)\tilde\eta(y)}\nonumber\\
&&\times\big(\delta((x\hspace{-0.5ex}+\hspace{-0.5ex}h)
-(y\hspace{-0.5ex}+\hspace{-0.5ex}h))-\delta((x\hspace{-0.5ex}+\hspace{-0.5ex}h)-y)
-\delta(x-(y\hspace{-0.5ex}+\hspace{-0.5ex}h))+\delta(x-y)\big)dxdy,
\end{eqnarray*}
where the spatial Dirac distributions come from the spatial white noise $\xi_{spat}$,
more precisely, they represent the covariance $\langle(\xi_{spat}^h-\xi_{spat})(x)
(\xi_{spat}^h-\xi_{spat})(y)\rangle_1$ of the increment $\xi_{spat}^h-\xi_{spat}$. 

\medskip

We integrate (\ref{l12}) against the weight $t+1$ in time over $t\in(-1,0)$.
This yields (\ref{l11}) once we establish the following three estimates:
The following estimate on the quadratic variation
\begin{equation}\label{l16}
\int\tilde\eta(\tilde\eta-\frac{1}{2}\tilde\eta^he^{-|h|}-\frac{1}{2}\tilde\eta^{-h}e^{-|h|})dx
\lesssim|h|,
\end{equation}
the following bound on the term under the time derivative
\begin{equation}\label{l15}
\int\tilde\eta(u^h-u)(1-\partial_x^2)^{-1}\tilde\eta(u^h-u)dx
\lesssim(e^{|h|}-1)^2\int\tilde\eta^2(u-c)^2dx,
\end{equation}
and the fact that ``elliptic term'' controls the desired term up to the term in (\ref{l15})
\begin{eqnarray}\label{l14}
\int\tilde\eta^2(u^h-u)^2dx&\le&
\frac{1}{C}\int\tilde\eta(u^h-u)(1-\partial_x^2)^{-1}\tilde\eta(-\partial_x^2)(\pi(u^h)-\pi(u))dx\nonumber\\
&+&C\int\tilde\eta(u^h-u)(1-\partial_x^2)^{-1}\tilde\eta(u^h-u)dx.
\end{eqnarray}

\medskip

We first address the quadratic variation term (\ref{l16}).
Writing
\begin{equation}\nonumber
\tilde\eta-\frac{1}{2}\tilde\eta^he^{-|h|}-\frac{1}{2}\tilde\eta^{-h}e^{-|h|}
=\tilde\eta(1-e^{-|h|})+e^{-|h|}(\tilde\eta-\frac{1}{2}\tilde\eta^h-\frac{1}{2}\tilde\eta^{-h})
\end{equation}
and performing a discrete integration by parts, we see that this term takes the form of
\begin{equation}\nonumber
(1-e^{-|h|})\int\tilde\eta^2dx+e^{-|h|}\int(\tilde\eta^h-\tilde\eta)^2dx,
\end{equation}
so that the estimate follows from the elementary estimate 
$\int(\tilde\eta^h-\tilde\eta)^2dx\le h^2\int(\partial_x\tilde\eta)^2dx$.

\medskip

We note that by duality, the estimate of the time-derivative term (\ref{l15}) is equivalent to
\begin{equation}\nonumber
\int\zeta\tilde\eta(u^h-u)dx
\lesssim(e^{|h|}-1)\left(\int(\zeta^2+(\partial_x\zeta)^2 dx
\int\tilde\eta^2(u-c)^2dx\right)^\frac{1}{2},
\end{equation}
which follows by discrete versions of integration by parts and Leibniz' rule
\begin{equation}\nonumber
\int\zeta\tilde\eta(u^h-u)dx=\int\big((\zeta^{-h}-\zeta)\tilde\eta
+\zeta^{-h}(\tilde\eta^{-h}-\tilde\eta)\big)(u-c)dx,
\end{equation}
Cauchy-Schwarz' inequality, the standard estimate
\begin{equation}
\int(\zeta^{-h}-\zeta)^2dx\le h^2\int(\partial_x\zeta)^2dx,
\end{equation}
and the following property of our cut-off function with exponential tails
\begin{eqnarray*}
\lefteqn{|\tilde\eta^{-h}(x)-\tilde\eta(x)|}\\
&=&\exp(-\frac{1}{2}\sqrt{x^2+1})
|\exp(\frac{1}{2}\sqrt{x^2+1}-\frac{1}{2}\sqrt{(x-h)^2+1})-1|\\
&\le&\tilde\eta(x)|\exp(\frac{|h|}{2})-1|.
\end{eqnarray*}

\medskip

Let us finally address the elliptic term (\ref{l14}). To this purpose we write (in our
physicist's way of omitting parentheses)
\begin{eqnarray*}
\lefteqn{\tilde\eta(-\partial_x^2)(\pi(u^h)-\pi(u))}\\
&=&(1-\partial_x^2)(\pi(u^h)-\pi(u))\tilde\eta+2\partial_x(\pi(u^h)-\pi(u))\partial_x\tilde\eta\\
&&-(\pi(u^h)-\pi(u))(1-\partial_x^2)\tilde\eta,
\end{eqnarray*}
so that by the symmetry of $(1-\partial_x^2)^{-1}$ (already used for (\ref{l12}))
\begin{eqnarray*}
\lefteqn{\int\tilde\eta(u^h-u)(1-\partial_x^2)^{-1}\tilde\eta(-\partial_x^2)(\pi(u^h)-\pi(u))dx}\\
&=&\int\tilde\eta^2(u^h-u)(\pi(u^h)-\pi(u))dx\\
&&-2\int(\partial_x\tilde\eta)(\pi(u^h)-\pi(u))\partial_x(1-\partial_x)^{-1}\tilde\eta(u^h-u)dx\\
&&-\int((1-\partial_x^2)\tilde\eta)(\pi(u^h))-\pi(u))(1-\partial_x^2)^{-1}\tilde\eta(u^h-u)dx.
\end{eqnarray*}
Using that the operators $\partial_x(1-\partial_x^2)^{-\frac{1}{2}}$ and $(1-\partial_x^2)^{-\frac{1}{2}}$
have operator norm 1 w.\ r.\ t.\ to $L^2$, we deduce the inequality (where we use the abbreviation
$\pi^h-\pi:=\pi(u^h)-\pi(u)$)
\begin{eqnarray}\nonumber
\lefteqn{\int\tilde\eta(u^h-u)(1-\partial_x^2)^{-1}\tilde\eta(-\partial_x^2)(\pi^h-\pi)dx}
\nonumber\\
&\ge&\int\tilde\eta^2(u^h-u)(\pi^h-\pi)dx\nonumber\\
&&-\Big(2\big(\int(\partial_x\tilde\eta)^2(\pi^h-\pi)^2dx\big)^\frac{1}{2}
+\big(\int((1-\partial_x^2)\tilde\eta)^2(\pi^h-\pi)^2dx\big)^\frac{1}{2}\Big)\nonumber\\
&&\times \left(\int\tilde\eta(u^h-u)(1-\partial_x^2)^{-1}\tilde\eta(u^h-u)dx\right)^\frac{1}{2}.
\end{eqnarray}
By the monotonicity properties (\ref{i.2}) of $\pi$, this yields
\begin{eqnarray}\label{l22}
\lefteqn{\int\tilde\eta(u^h-u)(1-\partial_x^2)^{-1}\tilde\eta(-\partial_x^2)(\pi(u^h)-\pi(u))dx}\nonumber\\
&\ge&\lambda\int\tilde\eta^2(u^h-u)^2dx\nonumber\\
&&-\Big(2\big(\int(\partial_x\tilde\eta)^2(u^h-u)^2dx\big)^\frac{1}{2}
+\big(\int((1-\partial_x^2)\tilde\eta)^2(u^h-u)^2dx\big)^\frac{1}{2}\Big)\nonumber\\
&&\times \left(\int\tilde\eta(u^h-u)(1-\partial_x^2)^{-1}\tilde\eta(u^h-u)dx\right)^\frac{1}{2}.
\end{eqnarray}
Our smoothing out of the exponential cut-off function ensures
\begin{equation}\label{l24}
|\partial_x\tilde\eta|+|\partial_x^2\tilde\eta|\lesssim\tilde\eta,
\end{equation}
which allows us to use Young's inequality in order to arrive at (\ref{l14}).

\bigskip

{\sc Proof of Lemma \ref{L2}}.
We will establish this lemma
in the strengthened version with 
the bulk average $\int_{-1}^0\int\eta udxdt$ replaced by the surface average
$c:=\int\eta udx_{|t=-1}$.
To this purpose we rewrite (\ref{i.1}) in form of
\begin{equation}\nonumber
\partial_t(u-c)=-(-\partial_x^2)(\pi(u)-\pi(c))+\xi.
\end{equation}
As in Lemma \ref{L1}, we replace $\eta$ by $\tilde\eta^2\sim \eta$ in the statement
of this lemma, with $\tilde\eta$ being 
the mollified version of $\eta_2$, cf.\ (\ref{l21}).
By the martingale argument we have like in Lemma \ref{L0}, cf.\ (\ref{l210}),
\begin{eqnarray*}
\lefteqn{\frac{d}{dt}\frac{1}{2}\Big\langle\int\tilde\eta(u-c)(1-\partial_x^2)^{-1}\tilde\eta(u-c)dx
\Big\rangle_1}\nonumber\\
&=&-\Big\langle\int\tilde\eta(u-c)(1-\partial_x^2)^{-1}\tilde\eta(-\partial_x^2)(\pi(u)-\pi(c))dx
\Big\rangle_1+\frac{1}{2}\int\frac{1}{2}\tilde\eta^2dx.
\end{eqnarray*}
%
We smuggle in an exponential term in the time variable with a rate $T\ll 1$
to be adjusted later:
\begin{eqnarray*}
\lefteqn{\frac{d}{dt}\exp(-\frac{t}{T})
\frac{1}{2}\Big\langle\int\tilde\eta(u-c)(1-\partial_x^2)^{-1}\tilde\eta(u-c)dx
\Big\rangle_1}\nonumber\\
&=&-\exp(-\frac{t}{T})\Big(\Big\langle
\frac{1}{2T}\int\tilde\eta(u-c)(1-\partial_x^2)^{-1}\tilde\eta(u-c)dx\nonumber\\
&&+\int\tilde\eta(u-c)(1-\partial_x^2)^{-1}\tilde\eta(-\partial_x^2)(\pi(u)-\pi(c))dx
\Big\rangle_1+\frac{1}{4}\int\tilde\eta^2dx\Big).
\end{eqnarray*}
Lemma \ref{L2} will follow from integration over $t\in(0,1)$ of this identity,
using the obvious estimates on the quadratic variation term
\begin{equation}\nonumber
\int\tilde\eta^2dx
\lesssim 1,
\end{equation}
and on the term under the time derivative
\begin{equation}\nonumber
\int\tilde\eta(u-c)(1-\partial_x^2)^{-1}\tilde\eta(u-c)dx
\le\int\tilde\eta^2(u-c)^2,
\end{equation}
once we show that the elliptic term controls the desired term for $T$ sufficiently small:
\begin{eqnarray}\label{l23}
\frac{1}{C}\int\tilde\eta^2(u-c)^2dx&\le&
\int\tilde\eta(u-c)(1-\partial_x^2)^{-1}\tilde\eta(-\partial_x^2)(\pi(u)-\pi(c))dx\nonumber\\
&+&\frac{1}{2T}\int\tilde\eta(u-c)(1-\partial_x^2)^{-1}\tilde\eta(u-c)dx.
\end{eqnarray}

\medskip

The argument for this estimate (\ref{l23}) 
on the elliptic term follows the lines of the one in Lemma \ref{L1}: Replacing the couple
$(u^h,u)$ from there by $(u,c)$, we arrive at
\begin{eqnarray}\nonumber
\lefteqn{\int\tilde\eta(u-c)(1-\partial_x^2)^{-1}\tilde\eta(-\partial_x^2)(\pi(u)-\pi(c))dx}\nonumber\\
&\ge&\lambda\int\tilde\eta^2(u-c)^2dx\nonumber\\
&&-\Big(2\big(\int(\partial_x\tilde\eta)^2(u-c)^2dx\big)^\frac{1}{2}
+\big(\int((1-\partial_x^2)\tilde\eta)^2(u-c)^2dx\big)^\frac{1}{2}\Big)\nonumber\\
&&\times \left(\int\tilde\eta(u-c)(1-\partial_x^2)^{-1}\tilde\eta(u-c)dx\right)^\frac{1}{2}.
\end{eqnarray}
Appealing to the estimates (\ref{l24}) of the smoothened exponential cut-off $\tilde\eta$
and Young's inequality, we obtain (\ref{l23}) for a sufficiently large $\frac{1}{T}$.

\bigskip

{\sc Proof of Lemma \ref{L4}}.
We fix an $r\le 1$ and note that the statement of this lemma follows from
\begin{eqnarray}\label{l410}
\lefteqn{\Big\langle\Big(\av_{-r^2}^0(\int\eta_ru
-\av_{-r^2}^0\int\eta_ru)^2dxdt\Big)^\frac{1}{2}\Big\rangle_r}\nonumber\\
&\lesssim&r^\frac{1}{2}
+\Big\langle\Big(\av_{-r^2}^0\int\eta_r(u-\int\eta_ru)^2dxdt\Big)^\frac{1}{2}\Big\rangle_r,
\end{eqnarray}
where $\langle\cdot\rangle_r$ denotes the expectation w.\ r.\ t.\ to the white noise restricted
to $(t,x)\in(-r^2,0)\times\mathbb{R}$, just by taking the expectation w.\ r.\ t.\ to $\langle\cdot\rangle_1$.
By the scale invariance (\ref{t1}) \& (\ref{t2}), it is thus sufficient to establish
the above for $r=1$. We shall replace the exponential averaging function
$\eta$ by its mollified version
\begin{equation}\nonumber
\tilde\eta(x)=\frac{1}{c_0}\exp(-\sqrt{x^2+1})\quad\mbox{with}\quad c_0:=\int\exp(-\sqrt{x^2+1})dx,
\end{equation}
noting that $\tilde\eta\sim\eta$
and pointing out the slight difference to Lemmas \ref{L1} and \ref{L2}, cf.\ (\ref{l21}). 
Indeed, $\tilde\eta\sim\eta$ is enough to replace $\eta$ by $\tilde\eta$ on the r.\ h.\ s.\
of (\ref{l410}); for the l.\ h.\ s.\ this follows the $L^2$-average in time of the estimate
\begin{eqnarray*}
\lefteqn{\Big|\int\eta udx-\int\tilde\eta udx\Big|}\\
&=&\Big|\int(\eta-\tilde\eta)(u-\int\eta u)dx\Big|
\lesssim\Big(\int\eta(u-\int\eta u)^2dx\Big)^\frac{1}{2}.
\end{eqnarray*}
Hence with the abbreviation $U(t):=\int\tilde\eta u dx$ we need to show that
\begin{eqnarray}\label{l44}
\Big\langle\Big(\int_{-1}^0(U
-\av_{-1}^0U)^2dxdt\Big)^\frac{1}{2}\Big\rangle_1
&\lesssim&1+
\Big\langle\Big(\int_{-1}^0\int\tilde\eta(u-U)^2dxdt\Big)^\frac{1}{2}\Big\rangle_1.
\end{eqnarray}

\medskip

After these preparations, we note that we may rewrite equation (\ref{i.1}) in form of
\begin{equation}\nonumber
\partial_tu=\partial_x^2\big(\pi(u)-\pi(U)\big)+\xi,
\end{equation}
From this we deduce the stochastic {\it ordinary} differential equation
\begin{equation}\nonumber
\partial_t\int\tilde\eta udx=\int\big(\pi(u)-\pi(U)\big)\partial_x^2\tilde\eta dx
+\sigma\partial_tW,
\end{equation}
where $W$ is a standard temporal Wiener process and the variance is given by
\begin{equation}\label{l43}
\sigma^2:=\int\tilde\eta dx\sim 1.
\end{equation}
We use the differential equation
in its time-integrated version
\begin{equation}\nonumber
\int_{-1}^0(\partial_t(U-\sigma W))^2dt=
\int_{-1}^0\Big(\int(\pi(u)-\pi(U))\partial_x^2\tilde\eta dx\Big)^2dt.
\end{equation}
Thanks to the Lipschitz continuity of $\pi$, (\ref{i.2}), and the fact that due
to our mollification $\tilde\eta$ of the exponential averaging function, we have
$|\partial_x^2\tilde\eta|\lesssim\tilde\eta$, this turns into the estimate
\begin{equation}\nonumber
\int_{-1}^0(\partial_t(U-\sigma W))^2dt\lesssim
\int_{-1}^0\int\tilde\eta(u-U)^2dxdt.
\end{equation}
By Poincar\'e's inequality (with vanishing mean value) and the triangle inequality,
and appealing to (\ref{l43}),
this turns into
\begin{equation}\nonumber
\int_{-1}^0(U-\int_{-1}^0U)^2dt\lesssim\int_{-1}^0(W-\int_{-1}^0W)^2dt
+\int_{-1}^0\int\tilde\eta(u-U)^2dxdt.
\end{equation}
By Jensen's inequality and the defining properties on the quadratic moments of the Brownian motion,
this implies
\begin{eqnarray*}
\lefteqn{\Big\langle\big(\int_{-1}^0(U-\int_{-1}^0U)^2dt\big)^\frac{1}{2}\Big\rangle_1}\nonumber\\
&\lesssim& \Big\langle\int_{-1}^0(W-\int_{-1}^0W)^2dt\Big\rangle^\frac{1}{2}
+\Big\langle\big(\int_{-1}^0\int\tilde\eta(u-U)^2dxdt\big)^\frac{1}{2}\Big\rangle_1\nonumber\\
&\lesssim& 1+\Big\langle\big(\int_{-1}^0\int\tilde\eta(u-U)^2dxdt\big)^\frac{1}{2}\Big\rangle_1,
\end{eqnarray*}
which is (\ref{l44}).

\bigskip

{\sc Proof of Lemma \ref{L3}}. 
%
First of all, the observable $\sup_{a_0 \in [\lambda,1]}\sup_{(t,x) \in (-1,0)\times\mathbb{R}}\eta(g^h-g)^2$ is 
independent from the noise $\xi$ outside of $(-1,0) \times \mathbb{R}$, 
so that we can replace the average $\langle \cdot \rangle_1$ by $\langle \cdot \rangle$.
\medskip

We start with the representation formula
\begin{equation}\label{l3.2}
g(a_0, t,x)=\int_{-1}^t\int   G(a_0(t-t'),x-x' ) \xi(t',x')dx'dt',
\end{equation}
where  $G(t,x):=\frac{1}{\sqrt{4\pi  t}}\exp(-\frac{|x|^2}{4  t})$ denotes the heat kernel for $\partial_t - \partial_x^2$. 
We simultaneously consider also 
\begin{align*}
 \partial_{a_0} g(a_0, t,x) = &\int_{-1}^t\int   (t-t') \partial_t  G(a_0(t-t'),x-x') 
  \xi(t',x')dx'dt',
\end{align*}
and 
%
%
argue that
\begin{align}
\langle(g(a_0,t,x)-g(a_0,s,y))^2\rangle\lesssim\sqrt{|t-s|}+|x-y|  \label{l3.1A}\\
\langle(\partial_{a_0} g(a_0,t,x)-\partial_{a_0} g(a_0,s,y))^2\rangle\lesssim\sqrt{|t-s|}+|x-y| \label{l3.1}
\end{align}
for all $(t,x),(s,y)\in\mathbb{R}\times\mathbb{R}$ and all $a_0 \in [\lambda,1]$.

\medskip
Because of the initial conditions and symmetry, we may w.\ l.\ o.\ g.\ assume that $-1\le s\le t$.
Using the defining property $\langle (\int \zeta \xi dxdt)^2 \rangle = \int \zeta^2 dxdt$ 
for a test function $\zeta$ of white noise
(that is, $\langle\zeta(t',x')\xi(s',y')\rangle=\delta(t'-s')\delta(x'-y')$ in the rough but efficient
physics language) we get an explicit expression for the covariances of these Gaussian fields. In the case of $g$ 
this is
\begin{eqnarray}
\langle g(a_0,t,x) g(a_0,s,y)\rangle
&=&\int_{-1}^s\int G(a_0(t-\tau),x-z)G(a_0(s-\tau),y-z)dzd\tau\notag\\
&=&\int_{-1}^s G(a_0(t+s-2\tau),x-y)d\tau \notag\\
&=&\frac{1}{2}\int_{t-s}^{t+s+2} G(a_0 \sigma,x-y)d\sigma, \label{Gaussian1}
\end{eqnarray}
where we used the semi-group property of $t\mapsto G(a_0 t,\cdot)$ in the middle identity.
In the of $\partial_{a_0}g$, i.e. the case where the integral kernel is given by $  \frac{d}{d a_0} G (a_0 t, x)$ 
the semi-group property has to be replaced by
\begin{align*}
 \int \frac{d}{da_0}  \frac{d}{da_0'} G (a_0 t, x) G(a_0' s, x-y) d x_{|a_0 = a_0'}  
&= \frac{d}{da_0}  \frac{d}{da_0'} G (a_0 t + a_0' s, y)_{|a_0 = a_0'}\\
&= t s \; \partial_t^2 G(a_0 (t+s),y).
\end{align*}
Using this formula, we obtain as in \eqref{Gaussian1}
\begin{align}
&\langle  \partial_{a_0} g(a_0,t,x)  \partial_{a_0}(a_0,s,y)\rangle\notag\\
&=\frac{1}{8}\int_{t-s}^{t+s+2} (\sigma^2-(t-s)^2)  \partial_t^2 G(a_0 \sigma,x-y)d\sigma.\notag
\end{align}
We now pass from covariance to increment. For  $g$ we get
\begin{eqnarray*}
\lefteqn{\langle(g(a_0,t,x)-g(a_0,s,y))^2\rangle}\nonumber\\
&=&\langle g^2(t,x)\rangle+\langle g^2(s,y)\rangle-2\langle g(t,x) g(s,y)\rangle\\
&=&\Big(\frac{1}{2}\int_{0}^{2(t+1)}+\frac{1}{2}\int_{0}^{2(s+1)}-\int_{t-s}^{(t+1)+(s+1)}\Big)G(a_0 \sigma,0)d\sigma\\
&&+\int_{t-s}^{t+s+2}(G(a_0\sigma,0)-G(a_0\sigma,x-y))d\sigma.
\end{eqnarray*}
By positivity and monotonicity of $G(a_0 \sigma,0)$ in $\sigma$ and $G(a_0\sigma,0)\ge G(a_0\sigma,z)$, this yields the inequality
\begin{eqnarray*}
\lefteqn{\langle(g(a_0,t,x)-g(a_0,s,y))^2\rangle}\nonumber\\
&\le&\int_0^{t-s}G(a_0\sigma,0)d\sigma
+\int_{0}^{\infty}(G(a_0\sigma,0)-G(a_0\sigma,x-y))d\sigma\\
&\le&\int_0^{t-s}G(a_0\sigma,0)d\sigma
+|x-y|\int_{0}^{\infty}(G(a_0\sigma,0)-G(a_0\sigma,1))d\sigma,
\end{eqnarray*}
where we used the scale invariance of $G(a_0\sigma,z)$ in the second step.
This inequality implies (\ref{l3.1A}) in this case  because of $G(a_0\sigma,0)\lesssim\sigma^{-\frac{1}{2}}$ and
$G(a_0\sigma,0)-G(a_0\sigma,1)\lesssim\min\{\sigma^{-\frac{1}{2}},\sigma^{-\frac{3}{2}}\}$.
Similarly we get for $\partial_{a_0}g$
\begin{eqnarray*}
\lefteqn{\langle(\partial_{a_0}g(a_0,t,x)-\partial_{a_0}g(a_0,s,y))^2\rangle}\nonumber\\
&=&\Big(\frac{1}{8}\int_{0}^{2(t+1)} +\frac{1}{8}\int_{0}^{2(s+1)} -\frac{1}{4}\int_{t-s}^{(t+1)+(s+1)}   \Big)\sigma^2 \partial_t^2 G(a_0 \sigma,0)d\sigma\\
&&+\frac{1}{4}\int_{t-s}^{(t+1)+(s+1)} (t-s)^2  \partial_t^2 G(a_0 \sigma,0)d\sigma \\
&&+\frac14 \int_{t-s}^{t+s+2} (\sigma^2-(t-s)^2 )  (\partial_t^2 G(a_0\sigma,0)-\partial_t^2 G(a_0\sigma,x-y))d\sigma\\
&\lesssim& \int_0^{t-s}\sigma^2 \sigma^{-\frac52} d\sigma + (t-s)^2 \int_{t-s}^{(t+1)+(s+1)} \sigma^{-\frac52} d \sigma \\
&&+ \int_0^\infty \sigma^2 (\partial_t^2 G(a_0\sigma,0)-\partial_t^2 G(a_0\sigma,x-y))d\sigma,
\end{eqnarray*}
so that the desired estimate \eqref{l3.1}  for $\partial_{a_0}g$ follows as well.

\medskip

We now apply Kolmogorov's continuity theorem to $g$ and $\partial_{a_0}g$; for the convenience of the reader
and because of its similarity to the proof of the main result of the paper,
we give a self-contained argument, to fix notation in the case of $g$. We first
appeal to Gaussianity to post-process (\ref{l3.1A}), which we rewrite as
\begin{equation}\nonumber
\langle\frac{1}{R}(g(t,x)-g(s,y))^2\rangle\lesssim 1\quad\mbox{provided}\;|t-s|\le 3R^2,|x-y|\le R
\end{equation}
for a given scale $R$.
We note that from (\ref{l3.2}) we see that the properties of being Gaussian and centered
transmits from $\xi$ to $\frac{1}{\sqrt{R}}(g(t,x)-g(s,y))$, so that by the above normalization
we have
\begin{equation}\label{l3.4}
\Big\langle\exp\big(\frac{1}{CR}(g(t,x)-g(s,y))^2\big)\Big\rangle\lesssim 1
\quad\mbox{for}\;|t-s|\le 3R^2,\;|x-y|\le R.
\end{equation}

\medskip

Our goal is to estimate exponential moments of the local H\"older-norm 
\begin{eqnarray*}\nonumber
[g]_{\alpha,(-1,0)\times(-1,1)}&:=&\sup_{R\in(0,1)}\frac{1}{R^\alpha}
\sup_{\stackrel{(t,x),(s,y)\in(-1,0)\times(-1,1)}{\sqrt{|t-s|}+|x-y|<R}}|g(t,x)-g(s,y)|,
\end{eqnarray*}
which amounts to exchange the expectation and the supremum over $(t,x)$, $(s,y)$ in (\ref{l3.4})
at the prize of a decreased H\"older exponent $\alpha<\frac{1}{2}$
To this purpose, we now argue that for $\alpha>0$, the supremum over a continuum can be replaced 
by the supremum over a discrete set: For $R<1$ we define the grid  
\begin{align*}
\Gamma_R= [-1,0]\times[-1,1]\cap(R^2\mathbb{Z}\times R\mathbb{Z})
\end{align*}
and claim that
\begin{eqnarray}
\lefteqn{[g]_{\alpha,(-1,0)\times(-1,1)}}\label{l3.7}\\
&\lesssim&\sup_{R}\frac{1}{R^\alpha}
\sup_{\stackrel{(t,x),(s,y)\in\Gamma_R }
{|t-s|\le 3 R^2,|x-y|\le R}}|g(t,x)-g(s,y)|=:\Lambda \nonumber,
\end{eqnarray}
where the first $sup$ runs over all $R$ of the form $2^{-N}$ for an integer $N \geq 1$.
Hence we have to show for arbitrary $(t,x),(s,y)\in(-1,0)\times(-1,1)$ that
\begin{equation}\label{l3.6}
|g(t,x)-g(s,y)|\lesssim\Lambda\big(\sqrt{|t-s|}+|x-y|\big)^\alpha.
\end{equation}
By density, we may assume that $(t,x),(s,y)\in r^2\mathbb{Z}\times r\mathbb{Z}$
for some dyadic $r=2^{-N}<1$ (this density argument requires the qualitative a priori
information of the continuity of $g$, which can be circumvented by approximating $\xi$). 
By symmetry and the triangle inequality, we may assume $s\le t$
and $x\le y$. For every dyadic level $n=N,N-1,\cdots $ we now recursively construct two sequences $(t_n,x_n)$
$(s_n,y_n)$ of space-time points (in fact, the space and time points can be constructed separately), 
starting from $(t_N,x_N)=(t,x)$ and
$(s_N,y_N)=(s,y)$, with the following properties 
\begin{enumerate}
\item[a)] they are in the corresponding lattice of scale $2^{-n}$, i.\ e.\
$(t_n,x_n),(s_n,x_n)\in (2^{-n})^2\mathbb{Z}\times 2^{-n}\mathbb{Z}$,
\item[b)] they are close to their predecessors in the sense of
$|t_{n}-t_{n+1}|,|s_{n}-s_{n+1}|\le 3(2^{-(n+1)})^2$ and $|x_{n}-x_{n+1}|,|y_{n}-y_{n+1}|\le 2^{-(n+1)}$,
so that by definition of $\Lambda$ we have
\begin{align}\label{l3.9}
|g(t_n,x_n)-g(t_{n+1},x_{n+1})|,& |g(s_n,y_n)-g(s_{n+1},y_{n+1})| 
\le \Lambda (2^{-(n+1)})^\alpha,
\end{align}
and 
\item[c)] such that $|t_n-s_n|$ and $|x_n-y_n|$ are minimized among these points.
\end{enumerate}

Because of the latter, we have
\begin{equation}\nonumber
(t_M,x_M)=(s_M,y_M)\quad\mbox{for some}\;M\;\mbox{with}\quad 2^{-M}\le\max\{\sqrt{|t-s|},|x-y|\},
\end{equation}
so that by the triangle inequality we gather from (\ref{l3.9})
\begin{equation}\nonumber
|g(t,x)-g(s,y)|\le\sum_{n=N-1}^M\Lambda (2^{-(n+1)})^\alpha\le\Lambda\frac{(2^{-M})^\alpha}{2^\alpha-1},
\end{equation}
which yields (\ref{l3.6}).

\medskip

Equipped with (\ref{l3.7}), we now may upgrade (\ref{l3.4}) to
\begin{equation}\label{l3.8}
\big\langle\exp\big(\frac{1}{C}[g]_{\alpha,(-1,0)\times(-1,1)}^2\big)\big\rangle\lesssim 1
\end{equation}
for $\alpha<\frac{1}{2}$. Indeed, (\ref{l3.7}) can be reformulated
on the level of characteristic functions as
\begin{equation}\nonumber
I\big([g]_{\alpha,(-1,0)\times(-1,1)}^2\ge M)
\le\sup_{R }\max_{(t,x),(s,y) \in \Gamma_R}I\big(\frac{1}{R}(g(t,x)-g(s,y))^2\ge\frac{M}{CR^{1-2\alpha}}\big),
\end{equation}
where as in (\ref{l3.7}) $R$ runs over all $2^{-N}$ for integers $N \geq 1$. Replacing the suprema by sums
in order to take the expectation, we obtain
\begin{equation}\nonumber
\big\langle I\big([g]_{\alpha,(-1,0)\times(-1,1)}^2\ge M)\big\rangle
\le\sum_{R}\sum_{(t,x),(s,y)}
\big\langle I\big(\frac{1}{R}(g(t,x)-g(s,y))^2\ge\frac{M}{CR^{1-2\alpha}}\big)\big\rangle.
\end{equation}
We now appeal to Chebyshev's inequality in order to make use of (\ref{l3.4}):
\begin{eqnarray*}\nonumber
\lefteqn{\big\langle I\big([g]_{\alpha,(-1,0)\times(-1,1)}^2\ge M)\big\rangle}\\
&\lesssim&\sum_{R}\sum_{(t,x),(s,y)}\exp(-\frac{M}{CR^{1-2\alpha}})\\
&\lesssim&\sum_{R}\frac{1}{R^3}\exp(-\frac{M}{CR^{1-2\alpha}})\\
&\stackrel{R\le 1,M\ge 1}{\le}
&\exp(-\frac{M}{C})\sum_{R}\frac{1}{R^3}\exp(-\frac{1}{C}(\frac{1}{R^{1-2\alpha}}-1))
\lesssim \exp(-\frac{M}{C}),
\end{eqnarray*}
where in the second step we have used that the number of pairs $(t,x),(s,y)$
of neighboring lattice points is bounded by $C\frac{1}{R^3}$ and in the last
step we have used that stretched exponential decay (recall $1-2\alpha>0$) beats polynomial growth.
The last estimate immediately yields (\ref{l3.8}).

\medskip

It remains to post-process (\ref{l3.8}) and the same bound for $\partial_{a_0}g$. We only need these bounds for second moments but
with the spatial origin replaced by any point $x$:
\begin{align}\label{l3.10A}
\big\langle[g]_{\alpha,(-1,0)\times(x-1,x+1)}^2\big\rangle +\big\langle[\partial_{a_0} g]_{\alpha,(-1,0)\times(x-1,x+1)}^2\big\rangle\lesssim 1
 \quad\mbox{for all}\;x\in\mathbb{R}.
\end{align}
We use these bounds in the embedding $H^{1} \hookrightarrow L^{\infty}$ to get 
\begin{align}
&\big\langle \sup_{a_0 \in [\lambda,1]} [g(a_0, \cdot, \cdot)]_{\alpha,(-1,0)\times(x-1,x+1)}^2\big\rangle\notag\\
 & \lesssim \Big\langle \int_{\lambda}^1  [g]_{\alpha,(-1,0)\times(x-1,x+1)}^2   + \int_{\lambda}^1[\partial_{a_0} g]_{\alpha,(-1,0)\times(x-1,x+1)}^2 \Big\rangle \notag\\
 & \lesssim 1
 \quad\mbox{for all}\;x\in\mathbb{R}. \label{l3.10}
\end{align}

\medskip
To obtain the statement of the lemma in form of
\begin{equation}\label{l3.11}
\big\langle \sup_{[\lambda,1]}\sup_{(-1,0)\times\mathbb{R}}\eta(g^h-g)^2\big\rangle\lesssim\min\{|h|^{2\alpha},1\}
\end{equation}
we distinguish the cases $|h|\le\frac{1}{2}$ and $|h|\ge\frac{1}{2}$.
In the first case we have for any $a_0$
\begin{eqnarray*}
\sup_{(-1,0)\times\mathbb{R}}\eta(g^h-g)^2
&\lesssim&
\sum_{x\in\mathbb{Z}}\exp(-|x|)\sup_{(-1,0)\times(x-\frac{1}{2},x+\frac{1}{2})}(g^h-g)^2\\
&\stackrel{|h|\le\frac{1}{2}}{\lesssim}&|h|^{2\alpha}\sum_{x\in\mathbb{Z}}\exp(-|x|)
[g]_{\alpha,(-1,0)\times(x-1,x+1)}^2,
\end{eqnarray*}
from which (\ref{l3.11}) follows by taking first the supremum over $a_0$ in each term and then the expectation and inserting
(\ref{l3.10}). In case of $|h|\ge\frac{1}{2}$, we proceed via 
\begin{eqnarray*}
\lefteqn{\sup_{(-1,0)\times\mathbb{R}}\eta(g^h-g)^2}\nonumber\\
&\lesssim&
\sup_{(-1,0)\times\mathbb{R}}\eta g^2+\sup_{(-1,0)\times\mathbb{R}}\eta^{-h} g^2\\
&\lesssim&
\sum_{x\in2\mathbb{Z}}(\exp(-|x|)+\exp(-|x-h|))\sup_{(-1,0)\times(x-1,x+1)}g^2\\
&\stackrel{g(t=-1)=0}{\lesssim}&
\sum_{x\in2\mathbb{Z}}(\exp(-|x|)+\exp(-|x-h|))\sup_{(-1,0)\times(x-1,x+1)}[g]_{\alpha,(-1,0)\times(x-1,x+1)}^2.
\end{eqnarray*}

\bigskip

{\sc Proof of Corollary~\ref{C1}}. We start by defining a modified local H\"older norm, based on the $D(u,r)$. For  $R >0$ set
\begin{equation}\label{co4}
\Gamma_R = \cap [-1,1]  \times [-1,1] \cap \big( R^2 \mathbb{Z} \times R \mathbb{Z} \big)  \,.
\end{equation}
Then we define the modified H\"older semi-norm
\begin{align*}
\llbracket u\rrbracket_\alpha = \sup_{R} \sup_{(\bar{t},\bar{x}) \in \Gamma_{R}} \frac{1}{\, R^{\alpha}}D(u^{(\bar{t},\bar{x})} ,R  )  ,
\end{align*}
where for a space-time point $(\bar{t},\bar{x})$ we write $u^{(\bar{t},\bar{x})}(t,x) = u(t+\bar{t},x+\bar{x})$ and the first supremum is taken over all $R = 2^{-N}$ for  integer $N \geq 1$. We claim that
\begin{equation}\label{co1}
[u]_\alpha \lesssim  \llbracket u\rrbracket_\alpha .
\end{equation}
This claim is established below, but first we proceed to prove Corollary~\ref{C1} assuming that \eqref{co1} holds.

\medskip

To this end, fix $\alpha < \alpha^{\prime} < \frac12$. From \eqref{t7} we get for any $R \in (0,1) $ and $1\leq \sigma< \infty$
\begin{align}
\Big\langle&  I \Big( \frac{1}{R^\alpha}D(u,R)   \geq \sigma \Big)  \Big\rangle \notag\\
& = \Big\langle  I \Big( \exp\Big( \Big( \frac{1}{R^{\alpha^{\prime}}}D(u,R) \Big)^{2 \frac{\an}{\alpha'}} \Big) \geq \exp\Big( ( \sigma R^{\alpha - \alpha^{\prime}} )^{ 2 \frac{\an}{\alpha'} } \Big) \Big)  \Big\rangle \notag\\
&\leq   \exp\Big( - \big( \sigma R^{\alpha - \alpha^{\prime}} \big)^{ 2 \frac{\an}{\alpha'} }  \Big)  \Big\langle   \exp\Big( \Big( \frac{1}{R^{\alpha^{\prime}}}D(u,R) \Big)^{2 \frac{\an}{\alpha'}} \Big)  \Big\rangle \notag\\
&\overset{\eqref{t7new}}{\lesssim}   \exp\Big( - \big( \sigma R^{\alpha - \alpha^{\prime}} \big)^{ 2 \frac{\an}{\alpha'} }  \Big)  \notag\\
&\lesssim   \exp\Big( - \frac{\sigma^{2 \frac{\an}{\alpha'}} }{C}  - \frac{1}{C}   R^{2 \frac{\an}{\alpha'}( \alpha - \alpha^{\prime} )}   \Big)  , \label{co3}
\end{align}
for a suitable constant $C$. In the third line we have used Chebyshev's inequality. By translation invariance the same bounds holds if $u$ is replaced by $u^{(\bar{t},\bar{x})}$ for any space-time point $(\bar{t},\bar{x})$. Therefore, we get
\begin{eqnarray*}
\Big\langle  I\big( \llbracket u\rrbracket_\alpha  \geq \sigma  \big) \Big\rangle &\leq &\sum_{R} \sum_{(\bar{t},\bar{x})\in \Lambda_{R}} \Big\langle I \Big( \frac{1}{R^\alpha} D(u^{(\bar{t},\bar{x})},R)   \geq \sigma  \Big) \Big\rangle \\
&\overset{\eqref{co3},\eqref{co4}}{\lesssim}  & \exp\Big( - \frac{\sigma^{2 \frac{\an}{\alpha'}} }{C} \Big)   \sum_{R} R^{-3}  \exp\Big( - \frac{1}{C}   R^{2 \frac{\an}{\alpha'}( \alpha - \alpha^{\prime} )}   \Big) \;\\
&\lesssim &  \exp\Big( - \frac{\sigma^{2 \frac{\an}{\alpha'}} }{C} \Big) .
\end{eqnarray*}
This fast decay of the tails of the distribution of the $ \llbracket u\rrbracket_\alpha$ implies the desired integrability property. 

\medskip

It remains to establish the bound \eqref{co1}. We rely on Campanato's characterization of H\"older spaces 
\cite[Theorem 5.5]{Giaquinta} which in our current context states that $[u]_\alpha$ is controlled by 
\begin{align}\label{cor13}
 \sup_{r <\frac12} \sup_{(t_0,x_0) \in [-1,1] \times [-1,1]} \frac{1}{\, r^{\alpha}}  \Big( \av_{-r^2}^0 \av_{-r}^r  \Big( u^{(t_0,x_0)} -  \av_{-r^2}^0 \av_{-r}^r  u^{(t_0,x_0)}  \Big)^2 \Big)^{\frac12} .
\end{align}
To see that $\llbracket u\rrbracket_\alpha $ controls this norm, we observe that for
$r>0$ satisfying $2^{-N-2} < r \leq 2^{-N-1}$ any arbitrary $(t_0,x_0) \in [-1,1] \times [-1,1]$ can be well approximated in $\Lambda_{2^{-N}}$, in the sense that there exists $(\bar{t},\bar{x}) \in \Lambda_{2^{-N}}$ satisfying $|x_0- \bar{x}|\leq 2^{-(N+1)}$ and  $|t_0-\bar{t}|\leq 2^{-2(N+1)}$. Then we get, for $R=2^{-N}$ using the definition of $\eta_{R}$
\begin{align*}
\frac{1}{\, r^{2\alpha}}&   \av_{-r^2}^0 \av_{-r}^r  \Big( u^{(t_0,x_0)} -  \av_{-r^2}^0 \av_{-r}^r  u^{(t_0,x_0)} dxdt \Big)^2 dxdt\\
 &\lesssim \frac{1}{R^{2\alpha}} \av_{-{R^2}}^0 \int \Big( \eta_{R} u^{(\bar{t},\bar{x}) } - \av_{-r^2}^0 \av_{-r}^r  u^{(t_0,x_0)}  dxdt \Big)^2 dx dt \\
&\lesssim  \frac{1}{R^{2\alpha}} D(u^{(\bar{t},\bar{x})},R)^2 \\
& \qquad + \frac{1}{R^{2\alpha}} \Big( \av_{-r^2}^{0} \av_{-r}^r  \Big( u^{(t_0,x_0)} - \av_{-R^2}^0 \int \eta_{R} u^{(\bar{t},\bar{x})} dx dt \Big) dxdt \Big)^2\\
&\lesssim  \frac{1}{R^{2\alpha}} D(u^{(\bar{t},\bar{x})},R)^2 .
\end{align*}

\medskip
Therefore, we can conclude that $\llbracket u\rrbracket_\alpha $ controls the Campanato norm defined in \eqref{cor13} and the proof of  \eqref{co1} (and therefore the proof of Corollary~\ref{C1}) is complete.


\begin{thebibliography}{99}

\bibitem{ArmstrongSmart}
S.~N. {Armstrong} and C.~K. {Smart}, Quantitative stochastic
  homogenization of convex integral functionals. ArXiv e-prints (2014).

\bibitem{AvellanedaLin}
M.~Avellaneda and F.-H. Lin, Compactness methods in the theory of homogenization. 
{\it Comm. Pure Appl. Math.} \textbf{40} (1987), no.~6, 803--847.

\bibitem{Conlon} J.\ Conlon, private communication (2015)

\bibitem{Debussche} A.\ Debussche, S.\ de Moor, M.\ Hofmanov\'a, A regularity result
for quasilinear stochastic partial differential equations of parabolic type.
{\it SIAM J. Math. Anal.}, {\bf 47} (2015), no.\ 2, 1590-1614 


\bibitem{FischerOtto}
J.~{Fischer} and F.~{Otto}, A higher-order large-scale regularity theory
  for random elliptic operators. ArXiv e-prints (2015).



\bibitem{FischerOtto2}
J.~{Fischer} and F.~{Otto}, Sublinear growth of the corrector in stochastic homogenization: 
Optimal stochastic estimates for slowly decaying correlations. ArXiv e-prints (2015).


\bibitem{Gess} B.\ Gess, Finite speed of propagation for stochastic porous media equations.
{\it SIAM J. Math. Anal.} {\bf 45} (2013), no. 5, 2734-2766




\bibitem{Giaquinta} M.\ Giaquinta, L.\ Martinazzi, An introduction to the regularity
theory for elliptic systems, harmonic maps and minimal graphs.
Appunti. Scuola Normale Superiore di Pisa (Nuova Serie), 
2. Edizioni della Normale, Pisa, 2005. xii+302 pp.

\bibitem{GloriaNeukammOtto}
A.~{Gloria}, S.~{Neukamm}, and F.~{Otto}, A regularity theory for random
  elliptic operators. ArXiv e-prints (2014).


\bibitem{LSU}
 O. A.\ Ladyzhenskaja, V. A.\ Solonnikov, N. N.\ Uraltseva, Linear and quasi-linear equations
of parabolic type,  New York, AMS, 1968.


\bibitem{Ledoux} M.\ Ledoux, Concentration of measure and logarithmic Sobolev inequalities. 
S\'eminaire de Probabilit\'es, XXXIII, 120–216, 
{\it Lecture Notes in Math.}, {\bf 1709}, Springer, Berlin, 1999



\bibitem{MarahrensOtto}
D.~{Marahrens} and F.~{Otto}, Annealed estimates on the Green
function, ArXiv e-prints (2013), online in {\it Probability Theory and Related Fields}. 

\bibitem{NaddafSpencer}
A.~Naddaf and T.~Spencer, On homogenization and scaling limit of some
  gradient perturbations of a massless free field. {\it Comm. Math. Phys.},
  {\bf 183} (1997), no.~1, 55--84.


\bibitem{Nash}
J.~Nash, Continuity of solutions of parabolic and elliptic equations.
{\it Amer. J. Math.}, {\bf 80} (1958), 931-954.

\end{thebibliography}
\end{document}